\documentclass[12pt]{amsart}
\usepackage[]{amsmath, amsthm, amsfonts, verbatim, amssymb, mathrsfs, mathtools, booktabs, multirow}

\usepackage[all]{xy}

\usepackage{hyperref}
\hypersetup{backref=true}
\usepackage{pifont}
\usepackage{enumerate}

\usepackage{graphicx}
\usepackage{epstopdf}
\usepackage{tikz,tikz-cd, color}
\usepackage[makeroom]{cancel}
\usetikzlibrary{arrows}
\usepackage{mathdots}
\usepackage{float}
\usepackage{caption, subcaption}
\usepackage{rotating}
\usepackage{adjustbox}
\usepackage{fullpage}

\def\Bl{{\rm Bl}}
\def\cD{{\mathcal D}}

\def\cO{{\mathcal O}}

\def\bP{\mathbb{P}}

\def\bZ{{\mathbb Z}}
\def\bQ{{\mathbb Q}}
\def\bC{{\mathbb C}}

\def\cA{{\mathcal A}}

\def\cC{{\mathcal C}}

\newtheorem {theo}{Theorem}
\newtheorem {coro}{Corollary}
\newtheorem {conj}{Conjecture}
\newtheorem {lemm}{Lemma}
\newtheorem {rem}{Remark}

\newtheorem {ques}{Question}
\newtheorem {prop}{Proposition}

\def\tM{\widetilde M}

\def\ts{\widetilde s}

\def\bs{\bigskip}
\def\ms{\medskip}
\def\ni{\noindent}

\def\pd{\partial}

\def\bQ{\mathbb{Q}}

\def\Pic{\rm{Pic}}
\def\bP{{\bf P}}

\def\oz1{d{\overline z}^1}
\def\oz2{d{\overline z}^2}
\def\oz3{d{\overline z}^3}

\def\ol{\overline l}

\def\oI{\overline I}

\def\oj{\overline j}

\def\oz{\overline z}

\def\oGamma{\overline \Gamma}

\def\oIq1{\oI_1\cdots\oI_{q-1}}
\def\oIq2{\oI_1\cdots\oI_{q-2}}

\def\dim{{\mbox{dim}}}

\def\Aut{{\rm Aut}}
\def\Hom{\mbox{Hom}}
\def\Ext{\mbox{Ext}}

\def\Sing{\rm Sing}
\def\NS{\rm NS}
\def\Ric{\rm Ric}
\newcommand{\HH}{\mathrm{H}}
\newcommand{\mult}{\rm mult}

\begin{document}

\author{Ching-Jui Lai}
\author{Sai-Kee Yeung}
\address{Department of Mathematics, National Cheng Kung University, Tainan 70101, Taiwan}
\address{Mathematics Department, Purdue University, West Lafayette, IN  47907 USA}              
\email{cjlai72@mail.ncku.edu.tw}
\email{yeungs@purdue.edu}

\title{Exceptional collection of objects on some fake projective planes}

\begin{abstract} The purpose of the article is to explain a new method to establish the existence of an exceptional collection of length three for a fake projective plane $M$ with non-trivial automorphism group, related to a conjecture of Galkin-Katzarkov-Mellit-Shinder in 2015. Our method shows that $30$ fake projective planes support such a sequence, most of which are new. In particular, this provides many new $H$-phantom categories. 
\end{abstract}

\maketitle

\tableofcontents


\section{Introduction}\label{Sec1}

\ni{\bf 1.1} A fake projective plane is a smooth compact complex surface $M$ with the same Betti 
numbers as $P_{\bC}^2$, but $M\ncong P_{\bC}^2$. This is a notion introduced by Mumford who also constructed 
the first example. All fake projective planes have recently been classified into $28$
classes by the work of Prasad-Yeung in \cite{PY}, where $60$ examples were constructed including a pair of examples 
for each class. Cartwright-Steger \cite{CS} confirmed a conjecture of \cite{PY} and proved that there were precisely $100$ fake projective planes from those $28$ classes, see also \cite{PY2}. It 
is known that a fake projective plane is a smooth complex two ball quotient, and has 
the smallest Euler number among smooth surfaces of general type.

Most of the fake projective planes have the property that the canonical line bundle 
$K_M$ can be written as $K_M=3L$, where $L$ is a generator of the N\'eron-Severi group, 
see Lemma \ref{cubicroot} for the complete list.  One motivation of the present article comes from
a question of Dolgachev and Prasad, who asked whether $\HH^0(M,2L)$ 
contains enough sections for geometric purposes, such as an embedding of $M$.

The other motivation comes from the recent research activities surrounding the search of 
exceptional collections and (quasi)-phantoms from the point of view of derived category, such as \cite{AO, BvBS, F, GS}  and \cite{GO}.
As for fake projective planes, this is related to whether $\HH^0(M,2L)$ is non-trivial, which has been  questioned and worked out in some examples in \cite{GKMS15}.  

\ms

\ni{\bf 1.2}  
Denote by $\cD^b(M)$ the bounded derived category of coherent sheaves on $M$. 
A sequence of objects $E_1, E_2, . . . , E_r$ of $\cD^b(M)$ is called an exceptional collection 
if $\Hom(E_j,E_i[k])$ is non-zero for $j\geq i$ and $k\in \bZ$ only when $i=j$ and $k=0,$ 
in which case it is one dimensional. In \cite{GKMS15} (or see \cite{GKMS}), the authors consider the problem of the existence 
of a special type exceptional collection on an $n$-dimensional fake projective space. 

\begin{conj}\label{GKMS} Assume that $M$ is an $n$-dimensional fake projective space with the canonical 
class divisible by $n+1$. Then for some choice of a line bundle $L$ such that $K_M=(n+1)L$, 
the sequence
$$\cO_M,-L,\dots,-nL$$
is an exceptional collection on $M$.
\end{conj}

In the cases of fake projective planes ($n=2$), it is easy to see that a necessary and sufficient condition 
for the above conjecture is to show that $\HH^0(M, 2L)=0$. This is proved in \cite{GKMS15} if $\Aut(M)$ has 
order 21. This is also proved for 2-adically uniformised fake projective planes in \cite{F}. 
The main result in this note aims to provide more examples to Conjecture \ref{GKMS}.

\newtheorem*{main}{Main Theorem}
\begin{main}\label{main} For $M$ a fake projective plane as listed in Table \ref{Tab:Main}, there is an $\Aut(X)$-invariant line bundle $L$ with $K_M=3L$ and the sequence $\mathcal{O}_M,-L,-2L$ forms an exceptional collection of $\cD^b(M)$. 

{\small
\centering
\begin{table}
   \captionof{table}{FPP with EC\label{Tab:Main}}
\begin{adjustbox}{center}$
\begin{array}{|c|c|c|c|c|c|}
\hline
\mbox{class}&M&{\rm Aut}(M)&\HH_1(M,\bZ)&H&\HH_1(M/H,\bZ)\\ \hline\hline
(a=7,p=2,\emptyset) &(a=7,p=2,\emptyset, D_3 2_7)&C_7:C_3&C_2^4&C_7&C_2\\ \hline
(a=7,p=2,\{7\})&(a=7,p=2,\{7\},D_3 2_7)&C_7:C_3&C_2^3&C_7&0\\ \hline
(\cC_{20},\{v_2\},\emptyset)&(\cC_{20},\{v_2\},\emptyset,D_3 2_7)&C_7:C_3&C_2^6&C_7&0\\  \hline
(\cC_2,p=2,\emptyset)&(\cC_2,p=2,\emptyset, d_3 D_3)&C_3\times C_3&C_2\times C_7&\Aut(M)&C_2\\  \hline
(\cC_2,p=2,\{3\})&(\cC_2,p=2,\{3\},d_3 D_3)&C_3\times C_3&C_7&\Aut(M)&0\\  \hline 
(\cC_{18},p=3,\emptyset)&(\cC_{18},p=3,\emptyset,d_3 D_3)&C_3\times C_3&C_2^2\times C_{13}&\Aut(M)&0\\  \hline
 (a=15,p=2,\{3,5\})&(a=15,p=2,\{3,5\},D_3)&C_3&C_3\times C_7&\Aut(M)&C_3\\ \cline{2-6}
 &(a=15,p=2,\{3,5\},3_3)&C_3&C_2^2\times C_3&\Aut(M)&C_3\\ \cline{2-6}
 &(a=15,p=2,\{3,5\},(D3)_3)&C_3&C_3&\Aut(M)&C_3\\ \hline    
(a=15,p=2,\{3\})&(a=15,p=2,\{3\},D_3)&C_3&C_2\times C_3\times C_7&\Aut(M)&C_2\times C_3\\ \cline{2-6}
&(a=15,p=2,\{3\},3_3)&C_3&C_2^3\times C_3&\Aut(M)&C_2\times C_3\\ \cline{2-6}
&(a=15,p=2,\{3\},(D3)_3)&C_3&C_2\times C_3&\Aut(M)&C_2\times C_3\\ \hline
(\cC_2,p=2,\emptyset)&(\cC_2,p=2,\emptyset, D_3 X_3)&C_3&C_2\times C_7\times C_9&\Aut(M)&C_2\times C_3\\ \cline{2-6}
&(\cC_2,p=2,\emptyset, (dD)_3 X_3)&C_3&C_2\times C_9&\Aut(M)&C_2\times C_3\\ \cline{2-6}
&(\cC_2,p=2,\emptyset, (d^2D)_3 X_3)&C_3&C_2\times C_9&\Aut(M)&C_2\times C_3\\ \hline      
\end{array}$
\end{adjustbox}
\end{table}
}\end{main}

There are 33 different pairs of fake projective planes with a non-trivial automorphism group. The above Table \ref{Tab:Main} covers $15$ pairs, while the other 18 pairs not covered by the Main Theorem are listed in Table \ref{Tab:unproven} of Section \ref{Sec9}, where we discuss the difficulty in our approach. As mentioned earlier, the results for the first three rows have been obtained earlier in \cite{GKMS15} by a different method. 
The Main Theorem is a combination of Theorem \ref{A} and Theorem \ref{B} to be explained in the next section. 

A consequence of Conjecture \ref{GKMS} is the existence of an $H$-phantom: A non-zero admissible subcategory $\cA$ of the derived category $\cD^b(M)$ is an $H$-phantom if the Hochschild homology ${\rm HH}_\bullet(\cA)=0$. 
From \cite[Corollary 1.2]{GKMS15} via taking $\cA$ to be the orthogonal complement in $\cD^b(M)$ of the exceptional sequence in the Main Theorem, we obtain 30 non-equivalent $H$-phantoms, cf. \cite{BO}.
\begin{coro} Any fake projective plane from the list of the Main Theorem admits an $H$-phantom in the derived categories $\cD^b(M)$.
\end{coro}
\begin{rem} 
For $M$ any fake projective plane as in the Main Theorem, take $G=C_3$, $C_7$, or $C_3:C_7$ to be a subgroup in $\Aut(M)$ and let $Z$ be the minimal resolution of $M/G$. By the same argument for proving \cite[Proposition 1.4, 1.6]{GKMS15}, we also know that both $\cD^b_G(M)$ and $\cD^b(Z)$ admit $H$-phantom subcategories. See \cite{GKMS15} for details.  
\end{rem}

\ms

\ni{\bf 1.3} To prove Conjecture \ref{GKMS}, the Riemann-Roch formula is not sufficient without an appropriate vanishing theorem, so the conjecture turns out to be rather subtle.

Our approach is geometric and different from \cite{GKMS15} and \cite{F}. We choose $L$ to be an $\Aut(M)$-invariant cubic root of $K_M$. The problem is reduced to a study of the geometry of invariant sections of $\HH^0(M,2L)$ if it exists. The proof relies on the classification of invariant curve and the group action on the fixed points on them. Our method depends mostly on the numerical property, and hence we propose the following slightly more general problem, which seems to be more accessible and still serves the purpose of searching for exceptional objects.

\begin{conj}\label{numGKMS} Assume that $M$ is an $n$-dimensional fake projective space with the canonical 
class {\bf numerically} divisible by $n+1$. Then for some choice of a line bundle $L$ such that $K_M=(n+1)L$ and a suitable choices of line bundles $E_i$'s with $E_i\equiv-iL$, $1\leq i\leq n$, the sequence $$\mathcal{O}_M,E_1,E_2,\dots, E_n$$ 
is an exceptional collection of $M$.
\end{conj}

 When $n=2$ and $\Aut(M)$ is large, namely, with order greater than $3$, by our method we can derive a contradiction to $h^0(M,2L)\neq0$. This holds in general and we find many exceptional collections. This implies immediately the following slightly stronger result (in the flavor of Conjecture \ref{numGKMS}). This follows from the proof of Theorem \ref{1} and \ref{2}, and is explained in the end of Section \ref{Sec4}.

\begin{theo}\label{A} Let $M$ be one of the fake projective plane in the list of the Main Theorem with $\Aut(M)=C_7:C_3$ or  $C_3\times C_3$. Suppose that $H=C_7$ when $\Aut(M)=C_7:C_3$ and $H=\Aut(M)$ otherwise. If $E_1$ and $E_2$ are two $H$-invariant torsion line bundles on $M$ and $L_i = L+E_i$, $i = 1,2$, then the sequence $\cO_M,-L_1,-2L_2$ forms an exceptional collection of $\cD^b(M)$.
\end{theo}

When $\Aut(M)=C_3$, we were not able to show directly that $h^0(M,2L)=0$. Instead, for the purpose of constructing exceptional objects, we assume that there are many invariant curves in the numerical class of $2L$ and derive a contradiction. This requires a study of the existence of two and three different invariant curves in the numerical class $2L$ and their possible intersection configurations. A careful analysis shows that there cannot be too many of them and leads to the required vanishing. In particular, we can prove the stronger Conjecture \ref{GKMS} when $M$ possesses enough $\Aut(M)$-invariant 3-torsions, cf. Corollary \ref{C3a} and \ref{C3b}.

\begin{theo}\label{B} Let $M$ be a fake projective plane with automorphism group $\Aut(M)=C_3$. If either 
\begin{enumerate}[(1)]
	\item $\HH_1(M/\Aut(M),\bZ)=C_3$, or 
	\item $\HH_1(M/\Aut(M),\bZ)=C_2\times C_3$ and $M$ is not in the class $\cC_{18}$,
\end{enumerate}
then for some $\Aut(M)$-invariant line bundle $L$ with $K_M=3L$, the sequence $\cO_M, -L, -2L$ forms an exceptional collection of $\cD^b(M)$.
\end{theo}
 
At this point, our argument is not sufficient to solve Conjecture \ref{numGKMS} for all fake projective planes with non-trivial automorphisms. For fake projective planes with non-trivial automorphism groups not covered in our theorems,  there are two classes:
\begin{enumerate}[$(a)$]
	\item $M$ is in class $\cC_{18}$, but $K_M\neq3L$ for any line bundle $L$. It is known that $K_M\equiv3H$ for some invariant line bundle $H$ and $M$ possesses many invariant 3-torsions.
	\item $M$ is not in class $\cC_{18}$, and $K_M=3L$ for a unique invariant line bundle $L$. Here $M$ possesses many nontrivial invariant $\Aut(M)$-torsions, but none of them has order 3.
\end{enumerate}

For class $\cC_{18}$, we are able to prove the following, cf. Section \ref{Sec8}.
\begin{prop}For a fake projective plane $M$ in the class $\cC_{18}$ with $\Aut(M)=C_3$, there is an invariant line bundle $L$ such that $K_M=3L+\omega$ for an invariant 3-torsion $\omega$ and $h^0(M,2L)=0.$ 
\end{prop}\
In particular, if furthermore $h^0(M,2L+\omega)=h^0(M,2L+2\omega)=0$, then the sequence $\cO_M,-L,-2L$ forms an exceptional collection by Lemma \ref{EC}. 
However, we are not able to show the existence of an exceptional collection as in Conjecture \ref{numGKMS} when either $h^0(M,2L+\omega)\neq0$ or $h^0(M,2L+2\omega)\neq0$. For fake projective planes in $(b)$, the difficulties encountered in our approach are explained in Section \ref{Sec9}.

After the completion of the first draft of the paper, the results in Theorem \ref{1} and \ref{2} of this paper (see Section \ref{Sec4} and \ref{Sec5}) were presented at the 4th South Kyushu Workshop on Algebra, Complex Ball Quotients and Related Topics, July 22-25, 2014, Kumamoto, Japan. The second author thanks Fumiharu Kato for his kind invitation. During the conference, J. Keum mentioned that he had obtained similar results to Theorem \ref{1} and \ref{2} as well, cf. \cite{K2}. 

\ms

\ni{\bf Organization.} This paper is organized as the following. In Section \ref{Sec2}, we study invariant line bundles on fake projective planes and the existence of invariant cubic root of the canonical class. In Section \ref{Sec3}, we prove the existence of invariant curves with fixed points. In Section \ref{Sec4} and \ref{Sec5}, we prove Theorem \ref{1}. In Section \ref{Sec6}, refining the study in Section \ref{Sec4} and \ref{Sec5}, we provide a list of possible configurations of two invariant curves in the numerical class $2L.$ In Section \ref{Sec7}, we show Conjecture \ref{GKMS} holds when $\Aut(M)=C_3$ and $\HH_1(M/\Aut(M),\bZ)=C_3$. In Section \ref{Sec8}, we study fake projective planes with $\Aut(M)=C_3$ and $\HH_1(M/C_3,\bZ)=C_2\times C_3$, and prove the remaining part of the Main Theorem. In the last Section \ref{Sec9}, we would explain the difficulties in applying our method to prove Conjecture \ref{numGKMS} for the remaining fake projective planes with non-trivial automorphisms. 

\ms

\ni{\bf Notation.} We work over $\bC$. Throughout this paper, we denote by $C_m$ the 
cyclic group of order $m$ and by $C_7:C_3$ the unique (up to isomorphism) nonabelian finite group  of order 21,
$$C_7:C_3=\langle x,y|x^3=y^3=1,xyx^{-1}=y^2\rangle.$$

The Picard group of a projective manifold $M$ is denoted by $\Pic(M)$, where we have $\sim$ the linear equivalence and $\equiv$ the numerical equivalence. The Neron-Severi group of $M$ is $\NS(M):=\Pic(X)/\equiv$ and $\NS(M)_\bQ:=\NS(M)\otimes_\bZ\bQ.$ We use additive notion for line bundles: $nL:=L^{\otimes n}$, and do not distinguish a line bundle $L$ with its associated Cartier divisor or $c_1(L).$ For two Cartier divisors $L_1, L_2$, we denote by $L_1\geqslant L_2$ if $L_1-L_2$ is an effective divisor. Also, we say $L_1$ is an $n$-th root of $L_2$ if $L_1=nL_2$ in $\Pic(X)$, and is a numerical $n$-the root of $L_2$ if $L_1\equiv nL_2$ in $\NS(X)$.

For a reduced proper curve $C$, we denote by $\nu:C^\nu\rightarrow C$ the normalization map. The sheaf $\delta:=\nu_*\cO_{C^\nu}/\cO_C$ is zero dimensional and supported on $\Sing(C)$. For the arithmetic genus $p_a(C):=h^1(C,\cO_C)$, we have $p_a(C)=g(C^\nu)+h^0(\delta)-s+1,$ where $s$ is the number of irreducible component of $C^\nu$ and $g(C^\nu)$ is the geometric genus.
  
\bs

\section{Line bundles on fake projective planes}\label{Sec2}

In this section, we study invariant line bundles on a fake projective plane and when does its canonical class admit an invariant cubic root. Recall that from \cite{PY} a fake projective plane is a ball quotient $M=B_\bC^2/\Pi$ for some lattice $\Pi\subseteq{\rm PU}(2,1)$, where $\Pi$ is constructed as a subgroup of a maximal arithmetic lattice $\oGamma\subset {\rm PU}(2,1)$ and $\Aut(M)=\oGamma/\Pi$.  We refer the reader to \cite{PY,CS} for details on the notations. 
The lattices $\oGamma$ and $\Pi$ are classified in \cite{PY,CS}

We remark that the Picard group $\Pic(M)=\NS(M)=\HH^2(M,\bZ)$ due to the cohomological
properties given in the definition of a fake projective plane.  We will use the fact throughout the following argument.

\ms

\ni{\bf 2.1}  Let $M$ be a fake projective plane.  First of all, we list all 
fake projective planes where $K_M$ has a cubic root as a line bundle.

\begin{lemm}\label{cubicroot}  Among the $100$ fake projective planes, $92$ of which satisfy the 
property that $K_M=3L$, where $L$ is a line bundle generating  $\NS(M)_\bQ.$
\end{lemm}

\ni{\bf Proof.} From the argument of \cite[\S10.2]{PY}, it is known that $K_M=3L$ if and only if $\Gamma$ can be lifted to become a lattice in ${\rm SU}(2,1)$, and $K_M=3L$ if $\HH^2(X,\bZ)$ has no 3-torsion. The latter fact is an immediate consequence of the Universal Coefficient Theorem, see {\bf 2.3} below or \cite[Lemma 3.4]{GKMS}. In \cite[\S10.2]{PY}, it also proves that $\Gamma$ can be 
lifted to ${\rm SU}(2,1)$ if the number field involved is not one of the types $\cC_{2}$ or $\cC_{18}$. 
There are $12$ candidates for $\Pi$ lying in $\cC_{2}$ or $\cC_{18}$. Out of these $12$ 
examples, $3$ of them do not have $3$-torsion elements in $\HH^2(M,\bZ)$ and hence the 
corresponding $\Pi$ can be lifted to ${\rm SU}(2,1)$. Finally, it is listed in the file \href{http://www.maths.usyd.edu.au/u/donaldc/fakeprojectiveplanes/registerofgps.txt}{registerofgps.txt} 
of the weblink of \cite{CS}, that the lattices can be lifted to ${\rm SU}(2,1)$ except for four cases in 
$\cC_{18}$, corresponding to $(\cC_{18},p=3,\{2\},D_3), (\cC_{18},p=3,\{2\},(dD)_3), 
(\cC_{18},p=3,\{2\},(d^2D)_3)$ and $(\cC_{18},p=3,\{2I\})$ in the notation of the file, 
see also Table 2 in \cite{Y2, Yadd}. Since there are two non-biholomorphic conjugate complex structures 
on such surfaces, it leads to the result that $92$ of the fake projective planes can be 
regarded as quotient of $B_{\bC}^2$ by a lattice in ${\rm SU}(2,1)$.
\qed

\ms 

\ni{\bf 2.2} Recall that for a fake projective plane $M$, the universal covering $\tM$ of $M$ is 
biholomorphic to $B_{\bC}^2$. Assume that $B_{\bC}^2$ is defined by a Hermitian form 
$F$ of signature $(2,1)$. Let ${\rm SU}(2,1)$ be the set of matrix elements in ${\rm GL}(3,\bC)$ 
preserving the Hermitian form $F$. Denote by $K_{\tM}$ the pull-back of the canonical 
line bundle on $M$ with respect to the universal covering map. Then $K_{\tM}$ is a 
${\rm SU}(2,1)$-equivariant holomorphic line bundle, and $K_{\tM}=3L_1$ in terms of 
a ${\rm SU}(2,1)$-equivariant holomorphic line bundle $L_1$ on $\tM$, cf. \cite{Ko}.

\begin{lemm}\label{invcubicroot} Let $M$ be a fake projective plane with $\Aut(M)\neq\{1\}$. 
\begin{enumerate}[$(a)$]
	\item Suppose that $M$ does not belong to the classes $\cC_2$ and $\cC_{18}$. Then $L_1$ 
descends as a holomorphic line bundle to $M$. Moreover, $L_1$ is invariant 
under $\Aut(M)$.  
	\item Suppose $M$ belongs to the class of $\cC_2$ or $\cC_{18}$ and is not one of the 
four cases of $\cC_{18}$ for which $\Pi$ cannot be lifted to ${\rm SU}(2,1)$. Then there 
is a subgroup $H<\Aut(M)$ of order $3$ for which $L_1$ is 
invariant under $H$.
\end{enumerate}
\end{lemm}
\ni{\bf Proof.} We begin with the proof of $(a)$. It is already proved in \cite{PY} that $\Pi$ can be lifted to ${\rm SU}(2,1)$, see Lemma 1 in \cite{CS}. 
From the set of generators of $\bar\Gamma$ listed by \cite{CS}, Cartwright and Steger actually show that $\bar\Gamma$ can 
be lifted to ${\rm SU}(2,1)$ as well. From Lemma \ref{cubicroot}, we already know that $L$ descend as a 
holomorphic line bundle to $M/\bar\Gamma$. Let $H$ be a subgroup of the automorphism 
group of $M$, then $M/H$ is a finite-sheeted covering of $M/\bar\Gamma$ from construction. 
Hence $L_1$ descends as a holomorphic line bundle from $\tM$ to $M/H$ as well, by 
pulling back from $\tM/\bar\Gamma$.

Consider now $\Pi$ belongs to the classes of $\cC_2$ or $\cC_{18}$ as in part $(b)$. 
From the file \href{http://www.maths.usyd.edu.au/u/donaldc/fakeprojectiveplanes/registerofgps.txt}{registerofgps.txt} in the weblink of \cite{CS}, we know that apart from the four cases of $\Pi$ in the table of Main Theorem, there is always a subgroup $H$ of the automorphism group of $M$ acting on $M$ such that the lattice associated to $M/H$  can be lifted to ${\rm SU}(2,1)$. Hence from the same argument as above, $L_1$ descends to a holomorphic 
line bundle to $M/H$. This implies that $L_1$ on $M$ is invariant under $H$.
\qed

\ms

In conclusion, for a fake projective plane $M$ with a non-trivial automorphism, if $M$ is not in class $\cC_{18}$, then $K_M=3L$ for an $\Aut(M)$-invariant line bundle $L$. 
Note that when $\Aut(M)=C_3\times C_3$, one can only find a cubic root $L$ of $K_X$ invariant under some $H=C_3<\Aut(M)$ from Lemma \ref{invcubicroot}. We will prove in Theorem \ref{2} 
that $L$ is indeed $\Aut(M)$-invariant. If $M$ is of classes $\cC_{18}$, then $\Aut(M)$ contains a subgroup $H=C_3$ from Table \ref{Tab:Main} and \ref{Tab:unproven}. 
By lifting a numerical cubic root of the canonical class of $M/H$, there still exists an $H$-invariant line bundle $L$ such that $K_M\equiv 3L$. 
There are two cases: 
\begin{enumerate}[(1)]
	\item $(\cC_{18},p=3,\emptyset, d_3,D_3)$: $\Aut(M)=C_3\times C_3$ and $K_M=3L$ for an $\Aut(M)$-invariant line bundle by Lemma \ref{cubicroot}, \ref{invcubicroot}, and proof of Theorem \ref{2}.
	\item $(\cC_{18},p=3,\{2\},D_3), (\cC_{18},p=3,\{2\},(dD)_3)$,  and $(\cC_{18},p=3,\{2\},(d^2D)_3)$: There is a unique $C_3$ factor in $\HH_1(M,\bZ)$, so $K_M$ indeed has three distinct numerical cubic roots, cf. Lemma \ref{invtor}. 
	But $K_M\neq 3L$ for any line bundle $L,$ cf. Lemma \ref{cubicroot}.
\end{enumerate}

\ms

We clearly have the following statement. 
\begin{lemm}\label{Autmodule} Suppose that $L$ is an $H$-invariant line bundle for some $H<\Aut(M)$. Then the space of sections $\HH^0(M,kL)$, if non-zero, is an $H$-module.  
\end{lemm}

\ms

\ni{\bf 2.3}\label{2.3} We consider torsion line bundles on $M$. The aim is to characterize when the canonical class $K_M$ of a fake projective plane  $M$ with $\Aut(M)=C_3$ has three distinct $\Aut(M)$-invariant cubic roots. 
This is crucial for our proof of the Main Theorem to be discussed from Section \ref{Sec6} to Section \ref{Sec9}. The key observation is that for such a surface $M$, there is a unique $C_3$ in its homology group $\HH_1(M,\bZ)$.

\begin{lemm}\label{invtor} Let $M$ be a fake projective plane admitting a nontrivial finite group $H=C_3<\Aut(M)$. 
If $\HH_1(M,\bZ)$ has exactly one copy of $C_3$ subgroup, then $\Pic(M)$ contains a subgroup $C_3$ consisting of $H$-invariant torsions. 
In particular, if $K_M=3L$ for some $L\in\Pic(M)$, then $K_M$ has three distinct cubic roots $L,L',L''$, which are 
$H$-invariant if so is $L$.
\end{lemm}
\ni{\bf Proof.} First we explain on a fake projective plane $M$, how torsion elements in $\HH_1(M,\bZ)$ corresponds to torsion elements in $\Pic(M).$

For a normal projective surface S, any holomorphic line bundle represents an 
element in $\NS(S)=i_*\HH^2(S,\bZ)\cap \HH^{1,1}(S)$, where $i:\bZ\rightarrow\bC$ 
is the inclusion map. In the case that $S$ is singular, we identify $\HH^{1,1}(S)$ with the corresponding part in $\HH^{1,1}(\widetilde{S})$
which is not contracted by $\mu$, where $\mu:{\widetilde{S}}\rightarrow S$ is the minimal resolution.
 Let us consider the torsion part of $\HH^2(S,\bZ)$. From the Universal Coefficient Theorem, 
we have the exact sequence 
$$0\rightarrow {\rm Ext}^1_{\bZ}(\HH_{1}(S,\bZ),\bZ)\rightarrow\HH^2(S,\bZ)\rightarrow \Hom_{\bZ}(\HH_2(S,\bZ),\bZ)\rightarrow0.$$
Since $\Hom_{\bZ}(\HH_2(S,\bZ),\bZ))$ is torsion free, 
for the sake of computation of torsion part of $i_*\HH^2(S,\bZ)\cap \HH^{1,1}(M)$, it suffices for us to investigate ${\rm Ext}^1_{\bZ}(\HH_{1}(S,\bZ),\bZ)$. 
On the other hand, for any abelian group $A$, we know that $\Ext_{\bZ}^1(\bZ/m\bZ,A)\cong A/mA.$ 
Hence $p$-torsions of $\HH^2(S,\bZ)$ corresponds to $p$-torsions of $\HH_{1}(S,\bZ)$.

The same argument applies to a fake projective plane $M$. For fake projective planes, all the torsion groups of $\HH_1(M,\bZ)$ are explicitly listed in the weblink associated to 
\cite{CS}.  The identification from the weblink together with the fact that $\Pic(M)\cong\HH^2(M,\bZ)$ conclude the proof of the first part.

We remark that for a fake projective plane $M$, the covering map $\pi:M\rightarrow S:=M/H$ is a Galois cover with isolated fixed points, cf. \cite{K1}. 
For a general smooth surface $M$ equipped with a finite automorphism group $H$ with isolated fixed points, 
there is a surjective group homomorphism 
$$\pi^*:\Pic(M/H)\rightarrow \Pic(M)^H=\{H-{\rm invariant\ line\ bundles}\}.$$
If we assume that $p$ is relative prime to the order of $H$,
an order $p$ real $1$-cycle on $S$ corresponds to an order $p$ real $1$-cycles on $M$ which is invariant under $H$.  In such case,
the pull-back of a non-trivial $p$-torsion line bundle from $S$ would still be non-trivial on $M$.

Now let $\langle\tau\rangle=\{0,\tau,2\tau\}\cong C_3<\Pic(M)$ be the subgroup of 3-torsions corresponding to the unique $C_3<\HH_1(M,\bZ)$. Here we use the additive notation on $\langle\tau\rangle$. 
If $g$ is a generator of $H=C_3<\Aut(M)$, then from our hypothesis $g\cdot\tau\in\langle\tau\rangle$. If $g\cdot\tau\neq \tau$, it has to be $2\tau$ or $0$.  
But $0$ is invariant under $\Aut(M)$ and hence $g\cdot\tau=2\tau$. As such $g\cdot(2\tau)=\tau$, and this implies that $g^2\cdot\tau=g\cdot(2\tau)=\tau$. 
But then $g^3\cdot\tau=2\tau\neq \tau$, a contradiction. In particular, $\{0,\tau,2\tau\}$ is a set of $H$-invariant torsion line bundles. 
If $K_M=3L$, then $L':=L+\tau$ and $L'':=L+2\tau$ are two other cubic roots of $K_M$. The rest is clear.
\qed

\bs

\section{Existence of invariant curves with fixed points}\label{Sec3}

In this section, $L$ is always a line bundle of $M$ such that $\NS(M)_\bQ=\left< L\right>.$ Note that $L^2=1$ by Poincar\'e duality. We also assume that the automorphism group $\Aut(M)$ of $M$ is non-trivial.

\ms

\ni{\bf 3.1}  We start with a simple statement, which has also been observed in \cite{GKMS15}. We include the proof for the convenience of the reader. 

\begin{lemm}\label{few} For a fake projective plane $M$, $h^0(M,2L)\leqslant 2.$
\end{lemm}

\ni{\bf Proof.}  Consider the homomorphism
$$\alpha:\HH^0(M,2L)\times \HH^0(M,2L)\rightarrow \HH^0(M,4L),$$
given by $\alpha(x,y)=x\times y$. This induces an injection
$$\mathbb{P}(\HH^0(M,2L))\times \mathbb{P}(\HH^0(M,2L))\stackrel{\beta}\rightarrow \mathbb{P}(\HH^0(M,4L)).$$
By \cite[Lemma 15.6.2]{Ko}, it follows that $h^0(M,4L)\geqslant 2h^0(M,2L)-1.$ Since $K_M\equiv 3L$ by the choice of $L$, $h^0(M,4L)=3$ by the Riemann-Roch formula and Kodaira vanishing theorem, and the lemma is proved.
\qed
\ms

For the induced action on $\HH^0(M,2L)$ when $L$ is invariant as in Lemma \ref{Autmodule}, the following key lemma proves the existence of an invariant curve equipped with a non-trivial group action when $h^0(M,2L)\neq0$. This is the cornerstone of our approach in this paper. 
\begin{lemm}\label{invcurve1} Let $M$ be a fake projective plane with $K_M\equiv3L$, where $L$ 
is invariant under a non-trivial cyclic subgroup $H<\Aut(M)$. If $h^0(M,2L)\neq0$, then there exists 
an $H$-invariant curve $\Sigma\sim 2L$ on which  $H$ acts non-trivially. Moreover, if $\Sigma$ is not irreducible and reduced, then one of the following holds: 
\begin{enumerate}[$(a)$]
	\item  $\Sigma=\Sigma_1+ \Sigma_2$, where $\Sigma_i\equiv L$ is irreducible and reduced for $i=1,2$. 
         In particular, $\Sigma_1$ and $\Sigma_2$ only intersect transversally at a smooth point. 
	\item  $\Sigma=2C$, where $C\equiv L$ is irreducible and reduced.
	\end{enumerate}
\end{lemm}
\ni{\bf Proof.} By Lemma \ref{few}, $h^0(M,2L)=1$ or 2. If $h^0(M,2L)=1$, then there exists a unique effective divisor $\Sigma\sim2L$. Since $h^*L=L$, we conclude that $h^*\Sigma=\Sigma$. Assume now $h^0(M,2L)=2$ so that there is an induced action of $H$ on $\mathbb{P}(\HH^0(M,2L))\cong P_{\bC}^1$. But the action of $H$ on $|2L|$ is linear and diagonalizable. Hence the existence of an invariant curve follows.

We claim that $H$ cannot act trivially on $\Sigma$. Assume on the contrary that it acts 
trivially on $\Sigma$. It follows that $\Sigma$ is fixed pointwise by $H$. Since $H$ is 
finite and $\Sigma$ is complex dimension $1$, we observe that $\Sigma$ must be totally 
geodesic. To see this, consider a real geodesic curve $c(t), |t|<\epsilon$ on $M$ with 
initial point $p\in \Sigma$ and initial tangent $\tau_p=c'(0)\in T_p\Sigma$. As both $p$ 
and $c'(0)$ are fixed by $H$, the whole geodesic curve $c(t), |t|<\epsilon$ is fixed by 
$H$ since the differential equation governing $c(t)$ is a second order ordinary equation 
and is determined by the initial conditions specified above. It follows that $c(t)$ actually 
lies on $\Sigma$. Since this is true for all points $p\in \Sigma$ and $\tau_p\in T_p\Sigma$, 
we conclude that $\Sigma$ is totally geodesic. On the other hand, from the result of \cite{PY}, we know that 
the lattice $\Pi$ associated to $M$ is arithmetic of second type.  It follows 
that there is no totally geodesic curve on $M$, cf. \cite[Lemma 8]{Y2}. The claim is proved.

Suppose that $\Sigma$ is not integral and write $\Sigma=\sum_im_i\Sigma_i$, where $\Sigma_i$'s 
are irreducible and reduced. Since $\Sigma_i\equiv n_iL$ for some $n_i\in\bZ_{>0}$ by $\NS(M)_\bQ=\left<L\right>$ and $\Sigma\equiv 2L$, we get $\sum_im_in_i=2$. Hence either $\Sigma=\Sigma_1+\Sigma_2$ with $\Sigma_i\equiv L$, 
or $\Sigma=2C$ with $C\equiv L$. Moreover, if $\Sigma=\Sigma_1+\Sigma_2$, then $\Sigma_1\cdot\Sigma_2=1$ 
and they can only intersect transversally at one smooth point. 
\qed

\ms
\ni{\bf 3.2} Now we apply holomorphic Lefschetz fixed point theorem to analyse the geometry of an $H$-invariant 
curve $\Sigma$ provided in Lemma \ref{invcurve1}. The main result is Proposition \ref{FP}, where we prove the existence of a fixed point. We will use the following lemma, cf. \cite{P}. 
\begin{lemm}\label{LFP} Let $C$ be a compact Riemann surface. Let $1\neq g\in{\rm Aut}(C)$ be an 
element of prime order $l$ acting non-trivially on $C$ with $n$ fixed points. Then for $\Delta=g(C)-\dim_\mathbb{C}\HH^{1}(\mathcal{O}_C)^{\rm inv}$, we have 
$$n=2-2g(C)+\frac{2l}{l-1}\Delta,$$
where $\HH^{1}(\mathcal{O}_C)^{\rm inv}$ is the eigenspace of eigenvalue $1$.         
\end{lemm}
\ni{\bf Proof.} We consider the holomorphic Lefschetz fixed point theorem,
\begin{align*}
\sum_{g\cdot p=p}\frac{1}{\det(1-\mathcal{J}_p(g^k))}={\rm tr}((g^k)^*|_{\mathrm{H}^0(C,\mathcal{O}_C)})-{\rm tr}((g^k)^*|_{\mathrm{H}^1(C,\mathcal{O}_C)}), \tag{$1$}
\end{align*}
where $\mathcal{J}_p(g^k)$ is the holomorphic Jacobian with respect to the action of $g^k$ at a fixed point $p$. 
We sum up $k=1,\dots,l-1$ of the above formula. 

Since $\HH^0(C,\cO_C)\cong\bC$, ${\rm tr}((g^k)^*|_{\mathrm{H}^0(C,\mathcal{O}_C)})=1$ for all $k.$ For the complex $\langle g\rangle$-module $V=\mathrm{H}^1(C,\mathcal{O}_C)$, since an eigenspace is one-dimensional, by considering the invariant and non-invariant part we deduce that 
$$\sum_{k=1}^{l-1}{\rm tr}((g^k)^*|_{\mathrm{H}^1(C,\mathcal{O}_C)})=(l-1)(g(C)-\Delta)-\Delta=(l-1)g(C)-l\Delta.$$ 
Hence the sum of the right hand sides of equation $(1)$ for $k=1,\dots, l-1$ equals to $l-1+l\Delta-(l-1)g(C)$.

For the left hand side of equation $(1)$, since $C$ is one-dimensional, $\mathcal{J}_p(g^k)=\rho^k$, 
where $\rho$ is an $l$-th root of unit. Hence each fixed point $p$ contributes  
$$\sum_{k=1}^{l-1}\frac{1}{1-\rho^k}=\frac{1}{2}(l-1),$$  
which then sums up to $\frac{n}{2}(l-1)$.  The equality in the lemma now follows easily.  

Here is an alternate argument, thanks to the suggestion of a referee.  Denote the quotient map by $\pi_g:C\rightarrow B$. Then from Serre duality $\HH^1(C,\cO_C)=\HH^0(C,K_C)^\vee$, we get
$$g(B)=h^0(B,K_B)=h^1(B,\cO_B)=\dim_\mathbb{C}\HH^{1}(\mathcal{O}_C)^{\rm inv}=g(C)-\Delta.$$ 
Now from the Riemann-Hurwitz formula, we get 
$$2g(C)-2=l\cdot(2g(B)-2)+\deg(R_{\pi_g})=l\cdot(2(g(C)-\Delta)-2)+n\cdot(l-1),$$
where $R_{\pi_g}$ is the ramification divisor. The lemma now follows.
\qed

\ms

We recall the following lemma, which is well-known to the experts.  
\begin{lemm}\label{Sch} For $C$ an irreducible and reduced curve on a fake projective plane $M$, $C$ is smooth of genus 3 if $C\equiv L$. If $C\equiv 2L$, then $g(C^\nu)\geq4$ and $h^0(\delta)\leq2.$
\end{lemm}
\ni{\bf Proof.} We first remark that for $C\subseteq M$, $g(C^\nu)\geq2$ as $M$ is hyperbolic. The  Ahlfors-Schwarz Lemma applied to the map induced by the normalization $\nu : C^\nu\rightarrow M$ cf. \cite{CCL}) for the manifolds equipped with Poincar\'e metrics implies that the K\"ahler forms satisfy $\nu^*\omega_M\leq\omega_{C^\nu}$, with equality if and only if it is a holomorphic isometry leading to totally geodesic $C$. Since there is no totally geodesic curve on a fake projective plane as mentioned in the proof of Lemma \ref{invcurve1}, the inequality is strict. Hence for $C\equiv kL$ with $k\geq1$, integrating over $C^\nu$, we get
$$2k= \frac{2}{3}(K_M\cdot C)<\deg(K_{C^\nu})=2g(C^\nu)-2=k(k+3)-2h^0(\delta),$$
where we used the fact that the Ricci curvature is $\frac{3}{2}$ of the holomorphic sectional 
curvature for the Poincar\'e metric on $M$ and the adjunction $p_a(C)=\frac{1}{2}C\cdot(K_M+C)$. 
Here in terms of the complex geodesic coordinates at the origin with $\frac{\pd}{\pd z_1}$ aligned with the tangential direction of $C$,
the Ricci curvature involved is 
$\Ric_{1\bar 1}=R_{1\bar 1 1\bar 1}+R_{1\bar 1 2\bar 2}$ and the holomorphic sectional curvature is
$R_{1\bar 1 1\bar 1}$, which equals $2R_{1\bar 1 2\bar 2}$ for the Poincar\'e metric on $B_{\bC}^2$,  with
the curvature tensor given by $R_{i\oj k\ol}=3\frac{\pd^4}{\pd z_i \pd \overline{z_j} \pd z_k \pd \overline{z_l} }\log (1-|z_1|^2-|z_2|^2)$,   cf. \cite{mok}.
Hence $k = 1$ implies that $h^0(\delta) = 0$ and $C$ is smooth with $g(C)=3$. The second statement is proved similarly.
\qed

 \begin{prop}\label{FP} Let $M$ be a fake projective plane with $K_M\equiv3L$. Suppose that $L$ is $H$-invariant for a non-trivial cyclic subgroup $H<\Aut(M)$ and $h^0(M,2L)\neq0$. Then there is an $H$-invariant curve $\Sigma\sim 2L$ with an $H$-fixed point. Moreover, one of the following holds:
	\begin{enumerate}[$(a)$]
		\item $\Sigma$ integral with $p_a(\Sigma)=6$, $g(\Sigma^\nu)\geq4$, and $h^0(\delta)\leq 2$;
		\item $\Sigma=2C$ and $C$ is smooth of genus 3;
		\item $\Sigma=\Sigma_1+ \Sigma_2$, where $\Sigma_i$'s are smooth of genus 3 and intersect transversally at a unique point.
	\end{enumerate}
\end{prop}

\ni{\bf Proof.} The existence of an $H$-invariant curve $\Sigma\sim 2L$ is from Lemma \ref{invcurve1}. Note that from \cite{PY}, $H$ can only be $C_3$ or $C_7$. To show the existence of an $H$-fixed point, we consider three cases: $\Sigma=\Sigma_1+\Sigma_2$, $\Sigma=2C$, or $\Sigma$ is irreducible and reduced as listed in Lemma \ref{invcurve1}.

If $\Sigma=\Sigma_1+\Sigma_2$, then $\Sigma_1\cap\Sigma_2=\{p\}$ is a point. As any element of $H$ carries an irreducible component of $\Sigma$ to another irreducible component and $|H|$ is odd, $\Sigma_i$'s are $H$-invariant and $p$ is an $H$-fixed point.

If $\Sigma=2C$, then $C\equiv L$ and is smooth of genus 3 by Lemma \ref{Sch}. 
If $H$ acts without fixed points on $C$, then the quotient $C/H$ is a compact Riemann surface of Euler-Poincar\'e number  
$$\chi_{\rm top}(C/H)=2-2g(C/|H|)=\frac{-4}{|H|}.$$ 
This is impossible for $|H|=3$ or $7$.

Suppose now that $\Sigma$ is irreducible and reduced. As proved in Lemma \ref{Sch}, $p_a(\Sigma)=6$, $g(\Sigma^\nu)\geq4$ and $h^0(\delta)\leq 2.$ If $h^0(\delta)\neq0$, then $0\neq|{\rm Sing}(\Sigma)|\leq 2$. Since the group action carries a singular point to a singular point and $|H|\geq3$ is odd, all the singular points are $H$-invariant. Suppose now $h^0(\delta)=0$ and hence $\Sigma$ is smooth. If $H$ acts without fixed points on $\Sigma$, then $\Sigma/H$ is a compact Riemann surface of Euler-Poincar\'e number  
$$\chi_{\rm top}(\Sigma/H)=2-2g(\Sigma/H)=\frac{-10}{|H|}.$$
This is impossible for $|H|=3$ or $7$.
\qed

\bs

\section{The case of  $\Aut(M)=C_7:C_3$}\label{Sec4}

In this section we prove the Main Theorem for a fake projective plane $M$ with ${\rm Aut}(M)=C_7:C_3$, which gives an alternate approach to such cases dealt with in \cite{GKMS15}.  
\begin{lemm}\label{ECeasy} Let $M$ be a fake projective with $K_M=3L$. The sequence $\mathcal{O}_M,-L,-2L$ forms an exceptional collection if and only if $h^0(M,2L)=0$.
\end{lemm}
\ni{\bf Proof.} This follows directly from the definition of an exceptional collection and the Serre duality, cf. \cite[Lemma 4.2]{GKMS15} or Lemma \ref{EC}.
\qed

\begin{theo}\label{1} Let $M$ be a fake projective plane with ${\rm Aut}(M)=C_7:C_3$. There is a line bundle $L$ such that $K_M=3L$ so that the sequence $\mathcal{O}_M,-L,-2L$ forms an exceptional collection.
\end{theo}
\ni{\bf Proof.} Let $L$ be any $\Aut(M)$-invariant cubic root of $K_M$ as given in Lemma \ref{invcubicroot}. By Lemma \ref{ECeasy}, we may assume that $\HH^0(M,2L)\neq\{0\}$.

Consider $H=C_7<{\rm Aut}(X)$, the unique $7$-Sylow subgroup. There is an $H$-invariant section $\Sigma\in\HH^0(M,2L)$ by Lemma \ref{invcurve1} and an $H$-fixed point by Proposition \ref{FP}. Moreover, $\Sigma$ is either irreducible and reduced, or $\Sigma=2C$, or $\Sigma=\Sigma_1+\Sigma_2$ is reducible. 

For the induced action of $H$ on $\Sigma^\nu$, observe that ${\rm Fix}(\Sigma)={\rm Fix}(M)\cap \Sigma$. In particular, $|{\rm Fix}(\Sigma)|\leq3$ by \cite{K1}. We denote \mbox{$x=\dim_\mathbb{C}\mathrm{H}^{1}(\mathcal{O}_{\Sigma^\nu})^{\rm inv}$} the dimension of $H$-invariant subspace and $n=|{\rm Fix}(\Sigma^\nu)|$ the number of $H$-fixed points on $\Sigma^\nu$. 

{\bf Case 1}: $\Sigma$ is irreducible and reduced. Here $p_a(\Sigma)=6$ and $h^0(\delta)\leq 2$ by Lemma \ref{Sch}.

Assume first that $\Sigma=\Sigma^\nu$, then $g(\Sigma)=6$ and $n\leq3$. For $l=7$, 
Lemma \ref{LFP} implies that $3n+7x=12$, where $(n,x)=(4,0)$ is the only nonnegative integer solution. This contradicts to the inequality $n\leq3$. 

Assume now that $\Sigma\neq\Sigma^\nu$. Applying Lemma \ref{LFP} to the lifted action of $H$ on $\Sigma^\nu$ with $l=|H|=7$, we get $3n+7x+h^0(\delta)=12$, where $h^0(\delta)=1$ or $2$.

If $h^0(\delta)=1$, then $3n+7x=11$ and there is no nonnegative integer solution. 

If $h^0(\delta)=2$, then $3n+7x=10$ and $(n,x)=(1,1)$. From the holomorphic Lefschetz fixed point theorem, 
we have  
$$ \frac{1}{1-\eta}+\xi_1+\xi_2+\xi_3=0,$$  
where $\eta, \xi_j\in(\bZ/7\bZ)^{\times}$. It can be checked directly from Matlab that there is no solution to the above equation.

{\bf Case 2}: $\Sigma=\Sigma_1\cup\Sigma_2$ is nodal at $p$ with two smooth irreducible components of genus 3. 

As observed in Lemma \ref{invcurve1}, $H$ acts on each $\Sigma_i$. Denote $n_i$ the number of $H$-fixed points on $\Sigma_i$. For $l=7$ in Lemma \ref{LFP}, we get  
$3n_i+7x_i=9$ and $(n_i,x_i)=(3,0)$ is the only solution in nonnegative integers. But then apart from the the fixed point $p$, there are two more fixed points on each $\Sigma_i$. This is a contradiction as $5=|{\rm Fix}(\Sigma)|>3$.

{\bf Case 3}: $\Sigma=2C$ with $C$ a smooth curve of genus 3.

Since $l=7$, $3n+7x=9$ by Lemma \ref{LFP}. It is only possible that $(n,x)=(3,0)$ and  
there are three smooth fixed points on $C$. From the holomorphic Lefschetz fixed point theorem, 
we have  
$$ \frac{1}{1-\eta_1}+\frac{1}{1-\eta_2}+\frac{1}{1-\eta_3}+\xi_1+\xi_2+\xi_3=1,$$  
where $\eta_i,\xi_j\in(\bZ/7\bZ)^{\times}$. It can be checked directly from Matlab 
that there is no solution to the above equation.
\qed

\ms
Theorem \ref{1} is a testing case of our geometric approach. However, the proof gets more complicated when the structure of $H$ gets simpler as we will see in later sections. 

\bs

\section{The case of  $\Aut(M)=C_3\times C_3$}\label{Sec5}

In this section, we prove the second part of the Main Theorem. From now on we assume 
$M$ is a fake projective plane with ${\rm Aut}(M)=C_3\times C_3$.  We refer the readers to the list in Section \ref{Sec1}. 

We have shown in Lemma \ref{cubicroot} and Lemma \ref{invcubicroot} that there is an ample line bundle $L$ with $K_M=3L$ 
and $L$ is $H$-invariant for some $H=C_3<{\rm Aut}(M)$. Hence $\HH^0(M,2L)$ is an ${\rm Aut}(M)$-module 
and we aim to show that $h^0(M,2L)=0$ as in Section \ref{Sec4}. 

The key point now is to refine our understanding of the singularities of an integral invariant curve $C\equiv 2L$ whenever it exists. Compare to Section \ref{Sec4}, the difficulty arises since the isotropic group of an invariant curve is now a smaller group $C_3$. Recall that from the result of \cite{CS} or \cite{K1}, a fixed point $o\in M$ of $H=C_3<\Aut(M)$ is of type $\frac{1}{3}(1,2)$. 
\begin{lemm}\label{local} Let $(C, o=(0,0))\subseteq\bC^2$ be an analytic germ of a singular reduced plane curve and $X_1,\dots, X_r$ be the irreducible branches of $C$ at $o\in C$. Then $$h^0(\delta)=\sum_ih^0(\delta_i)+\sum_{i<j}(X_i.X_j).$$ 
In particular, $h^0(\delta)\geq r(r-1)/2$. 
	
Furthermore, suppose that $H=C_3$ acts on $\bC^2$ with weight $\frac{1}{3}(1,2)$ and $(C, o)$ is $H$-invariant. If the induced action on $(C,o)$ is nontrivial and $h^0(\delta)\leq2$, then either 
	\begin{enumerate}
		\item[$(a)$] $h^0(\delta)=1$ and $o\in C$ is a node, or
		\item[$(b)$] $h^0(\delta)=2$ and $o\in C$ is a tacnode.
	\end{enumerate}
In particular, in both cases $r=2$ and $o\in C$ lifts to two $H$-fixed points on $C^\nu$.	
\end{lemm}
\ni{\bf Proof.} First part is given in Hironaka \cite{Hi}. For the second part, we first observe that $h^0(\delta)\leq 2$ implies that $r\leq2$.

Suppose that $r=1$. We consider the sequence 
$$0\rightarrow\cO_{C,o}\cong\frac{\bC[[x,y]]}{(f(x,y))}\xrightarrow{\phi} \bC[[t]]\rightarrow\delta\rightarrow0,$$
where $f(x,y)$ is the defining equation of $C$, $\phi(x)=u(t)=\sum_{m\geq0}u_mt^m$, and $\phi(y)=v(t)=\sum_{n\geq0}v_nt^n$. 
Here we choose $(x,y)$ to be $H$-invariant coordinate with $\omega\cdot x=\omega x$ and $\omega\cdot y=\omega^2y$, where $\omega=\exp(2\pi i/3)$. 
Up to an analytic change of coordinates we can assume that $\omega\cdot t=\omega^\alpha t$ for $\alpha\in\{1,2\}$. Since $\phi$ is an $H$-invariant $\bC$-algebra homomorphism, we have 
\begin{align*}\begin{cases}\omega u(t)=\phi(\omega x)=\phi(\omega\cdot x)=\omega\cdot u(t)=\sum_{m\geq0}u_m\omega^{\alpha m}t^m\\
\omega^2 v(t)=\phi(\omega^2y)=\phi(\omega\cdot y)=\omega\cdot v(t)=\sum_{n\geq0}v_n\omega^{\alpha n}t^n\end{cases}.
\end{align*}
Hence for any nonzero $u_m$ and $v_n$, we have 
$$\alpha m\equiv 1,\ \alpha n\equiv 2\ \mod 3.$$ 

Assume that $\alpha=1$ and write 
\begin{align*}(u(t),v(t))=\left(t(\sum_{m\geq0}a_mt^{3m}), t^2(\sum_{n\geq0}b_nt^{3n})\right).
\end{align*}
Here $a_0=0$ or otherwise $C$ is smooth. Furthermore, if $b_0\neq0$, then $v(t)=t^2\cdot{\rm unit}$. Hence the $k$-algebra $\cO_{C,o}$ is of the form $k[[t^2\cdot{\rm unit}, t^{4+3l}\cdot{\rm unit}]]$ for some $l\geq0$. 
But then $\delta$ contains at least $t, t^3, t^5$, which contradicts to $h^0(\delta)\leq2$.  If $b_0=0$, then the same computation leads to $h^0(\delta)>2$, which is a contradiction. The case $\alpha=2$ is similar. Hence we must have $r>1.$

Suppose now that $r=2$ so that $h^0(\delta)=h^0(\delta_1)+h^0(\delta_2)+X_1\cdot X_2\geq X_1\cdot X_2\geq1$. 
If $h^0(\delta)=1$, then both $X_i$'s are smooth and intersect transversally. This is the nodal case $(a)$. 
Assume that $h^0(\delta)=2$. If $h^0(\delta_1)=h^0(\delta_2)=0$ and $X_1\cdot X_2=2$, then $X_1=\{x=0\}$ after a change of coordinates and $X_2\cong \{x-y^2=0\}$. We get $C\cong\{x(x-y^2)=0\},$ which is case $(b)$. On the other hand, if $h^0(\delta_1)+h^0(\delta_2)=1$ and assume that $X_1$ is smooth, then $X_1\cdot X_2\geq2$ as $X_2$ is singular. This is impossible. 
\qed

\ms

We now refine Lemma \ref{FP} in the case when an invariant curve $C$ is irreducible and reduced. 
\begin{lemm}\label{invirr} Let $M$ be a fake projective plane with $K_M\equiv3L$ 
where $L$ is $H$-invariant for some $H=C_3<\Aut(M)$. Let $C\equiv kL$ be an integral $H$-invariant curve with $k=1$ or $2$. Denote by $n$ the number of $H$-fixed points on $C^\nu$ and $x=\dim_\bC\HH^1(C^\nu,\cO)^{\rm inv}$. There is a 
finite list of $C$ according to the triple $(n,h^0(\delta),x)$: 
\begin{enumerate}
	\item[$(N)$] $(n,h^0(\delta),x)=(2,0,1)$: $C\equiv L$ is smooth with two smooth fixed points;
	\item[$(I_1)$] $(n,h^0(\delta),x)=(2,0,2)$: $C\equiv 2L$ is smooth of $g(C)=6$ with two smooth fixed points;
	\item[$(I_2)$] $(n,h^0(\delta),x)=(4,1,1)$: $C\equiv 2L$ has one fixed node, which is the unique singularity of $C$, and two smooth fixed points;
	\item[$(I_3)$] $(n,h^0(\delta),x)=(3,2,1)$: $C\equiv 2L$ has one fixed tacnode and one smooth fixed point. 	         
\end{enumerate}
\end{lemm}

\ni{\bf Proof.}  From Lemma \ref{LFP} and \ref{Sch}, we have either
\begin{enumerate}[(1)]
	\item $C\equiv L$: $h^0(\delta)=0$, $g(C)=3$, and $n+3x=5$, or 
	\item $C\equiv 2L$: $h^0(\delta)\leq2$, $g(C^\nu)\geq4$, and $n+h^0(\delta)+3x=8$.
\end{enumerate}
Here $n\geq1$ by Proposition \ref{FP}. Observe that from the proof of Proposition \ref{FP}, all the singularities of $C$ are $H$-fixed points. Note that the set of fixed points of $C$ satisfies $|{\rm Fix}(C)|=|{\rm Fix}(M)\cap C|\leq|{\rm Fix}(M)|=3$ by the work of \cite{K1}. 

In case $(1)$, since $C=C^\nu\subseteq M$ and $n=|{\rm Fix}(C)|\leq3$, 
there is only one solution $(n,x)=(2,1)$.\footnote{Note that there is no contradiction to 
holomorphic Lefschetz fixed point theorem as 
$$\frac{1}{1-\omega}+\frac{1}{1-\omega^2}+\omega+\omega^2=0,$$
where $\omega=\exp(\frac{2\pi i}{3})$.} This is the case $(N)$.

In case (2), we have the following possible solutions:
{\large
$$\begin{array}{|c|c|c|}
\hline
\delta & {\rm equation} & (n,x)  \\ \hline
h^0(\delta)=0 & n+3x=8 & (8,0), (5,1), (2,2)  \\ \hline
h^0(\delta)=1 & n+3x=7 & (7,0), (4,1), (1,2)  \\ \hline
h^0(\delta)=2 & n+3x=6 & (6,0), (3,1) \\ \hline
\end{array}$$
}
If $h^0(\delta)=0$, then $C=C^\nu$ and $n=|{\rm Fix}(C)|\leq3$. Hence $(n,x)=(2,2)$ and this is case $(I_1)$.

By Lemma \ref{local}, $h^0(\delta)=1$ occurs only when the unique singular point is a node, which then lifts to two $H$-fixed points on $C^v$. Hence $n\geq2$ and $(4,1)$ is the only solution as there are at most two more smooth fixed points by $|{\rm Fix}(C)|\leq3$. This gives case $(I_2)$. 

If now $h^0(\delta)=2$, then $|{\rm Sing}(C)|=1$ or 2. If there are two singular points, then by Lemma \ref{local} these are two $H$-fixed nodes. These two nodes lift to four $H$-fixed points on $C^\nu$ and $n\geq4$. Hence $(n,x)=(6,0)$ is the only solution. But then there must be two more smooth $H$-fixed points on $C$ and this contradicts to $|{\rm Fix}(C)|\leq3.$ 

If there is only one singular point, then by Lemma \ref{local} it again lifts to two $H$-fixed points on $C^\nu$ and $n\geq 2$. If $(n,x)=(6,0)$, then there must be four more 
smooth $H$-fixed points on $C$, which contradicts to $|{\rm Fix}(C)|\leq3$. Hence $(n,x)=(3,1)$ and this is case $(I_3)$.\footnote{Note that the holomorphic Lefschetz fixed point theorem has a solution,
$$\frac{1}{1-\omega}+\frac{1}{1-\omega}+\frac{1}{1-\omega}+\omega+\omega^2+\omega^2=0.$$} 
\qed
\ms

We make the following simple observation.
\begin{lemm}\label{faith} There is no non-trivial faithful action of $C_3\times C_3$ fixing a point on any fake projective plane.
\end{lemm}
\ni{\bf Proof.} From the work of \cite{PY} and \cite{CS}, the automorphism group of $M=B_{\bC}^2/\Pi$ is given
by the quotient group $H=\bar\Gamma/\Pi$. In the cases that  $H$ contains $C_3\times C_3$, actually $H=C_3\times C_3$.  It is shown case 
by case in the files of the weblink associated to \cite{CS} that the singularities of $H$ consists of $12$ points on $M$, each being a fixed point
of one of the four $C_3$ subgroups of $C_3\times C_3$.  In particular, there is no point on $M$ fixed by all elements of $H$.

Alternatively, we observed that the finite group does not contain any subgroup acting as complex reflections on a fake projective plane $M$: a fixed curve of a complex reflection is totally geodesic, which does not exist on a fake projective plane, cf. Lemma \ref{invcurve1}.
The action is an ${\rm SL}(2,\bC)$ action since it preserves the K\"ahler-Einstein volume form.  We may then resort to the classification of actions of finite subgroup of ${\rm SL}(2,\bC)$ on $\bC^2$ as given in \cite[Corollary 4-6-16]{ma} to conclude the proof.
\qed

\begin{theo}\label{2} Let $M$ be a fake projective plane with $\Aut(M)=C_3\times C_3$. 
There is an $\Aut(M)$-invariant line bundle $L$  such that $K_M=3L$, and the sequence 
$\cO_M,-L,-2L$ forms an exceptional collection.
\end{theo}
\ni{\bf Proof.} From the classification of fake projective plane, a fake projective plane $M$ with $\Aut(M)=C_3\times C_3$ are all listed in the table of the Main Theorem and satisfies Lemma \ref{invcubicroot} $(2)$. 
Hence, there is an $H$-invariant line bundle $L$ such that $K_M=3L$ for some $H=C_3<\Aut(M)$. We prove that this $L$ is indeed $\Aut(M)$-invariant:  
For $g\in\Aut(M)$, since $K_M=3L$ and 
$$3(g\cdot L)=g\cdot (3L)=g\cdot K_M=K_M=3L,$$ 
we see that $g\cdot L-L$ is a 3-torson line bundle. But from the proof of Lemma 4, a torsion line bundles of $M$ corresponds to a torsion elements in $\HH_1(M,\bZ)$, which as we can find 
in the table of the Main Theorem that its order can never be 3. Hence $L$ is $\Aut(M)$-invariant. 

Suppose that $H^0(M,2L)\neq0$ and let $\Sigma$ be an $\Aut(M)$-invariant section. 
Note that $\Aut(M)$ has four $C_3$ subgroups, denoted by $G_1,\dots,G_4$. From the proof of Lemma \ref{faith} (or cf. \cite{K1}), there are twelve points $P_i, i=1,\dots,12,$ of $M$ and each point is fixed by some $G_i$. The stabilizer 
of each $P_i$ is $G_j$ for some $j$. Hence on the quotient surface $Y:=M/\Aut(M)$, there are four points 
of type $\frac13(1,2)$.  

Let $G_1$ be the first $C_3$ factor and $G_2$ be the second $C_3$ factor. 

Consider $G_1$-action on $\Sigma$. From  Proposition \ref{FP} and Lemma \ref{invirr}, there are three possibilities:
\begin{enumerate}
	\item $\Sigma$ is integral and the number of smooth fixed points is at most two. 
	\item $\Sigma=2C$ and $C$ is smooth of genus 3 with two smooth $G_1$-fixed points.
	\item $\Sigma$ is reduced with two smooth components $\Sigma_1$ and $\Sigma_2$ of genus 3. Moreover, $G_1$ acts on each component $\Sigma_i$ with two smooth fixed points. 
\end{enumerate}

Since $\Sigma$ is $\Aut(M)$-invariant, the curve $\Sigma$ in (1), $C$ in (2), and $\Sigma_i$ in (3) are all invariant under $G_2$. Moreover, $G_2$-action permutes $G_1$-fixed points by Lemma \ref{faith}. 
Since each curve $\Sigma$, $C$, or $\Sigma_i$ has at least 1 smooth $G_1$-fixed point $P$, 
there are at least three $G_1$-fixed points as the $G_2$-orbit of $P$ on them. This is a contradiction to the above list of possible $\Sigma$. 
\qed

\ms
Theorem \ref{A} is the combination of Theorem \ref{1} and Theorem \ref{2}. 
\ms

\ni{\bf Proof.} (of Theorem \ref{A}) For $\cO_M, -L_1,-2L_2$ in Theorem \ref{A} to form an exceptional collection, we need to show that
$$h^i(M,L_1)=h^i(M,2L_2)=h^i(M,2L_2-L_1)=0,\ i=0,1,2.$$
We consider vanishing of $h^i(M, L_1)$ first. Note that $h^2(M, L_1) = h^0(M, K_M-L_1).$ 
Since both $L_1$ and $K_M-L_1$ are invariant under $\Aut(M)$, $h^0(M,2L_1)=0=h^0(M,K_M-L_1)$ 
by the same proof as in Theorem \ref{1} and \ref{2}. It follows that $h^0(M,L_1)=0$ and then by the Riemann-Roch formula $h^1(M,L_1) = 0$. The other vanishing are proved similarly. 
\qed

\bs

\section{Invariant curves on $M$ when $\Aut(M)=C_3$}\label{Sec6}

Throughout this section, we assume that $M$ is a fake projective plane with $\Aut(M)=C_3$ unless otherwise stated.

\ms

\begin{figure}
\centering
\caption{INVARIANT CURVES}
\includegraphics*[scale=1]{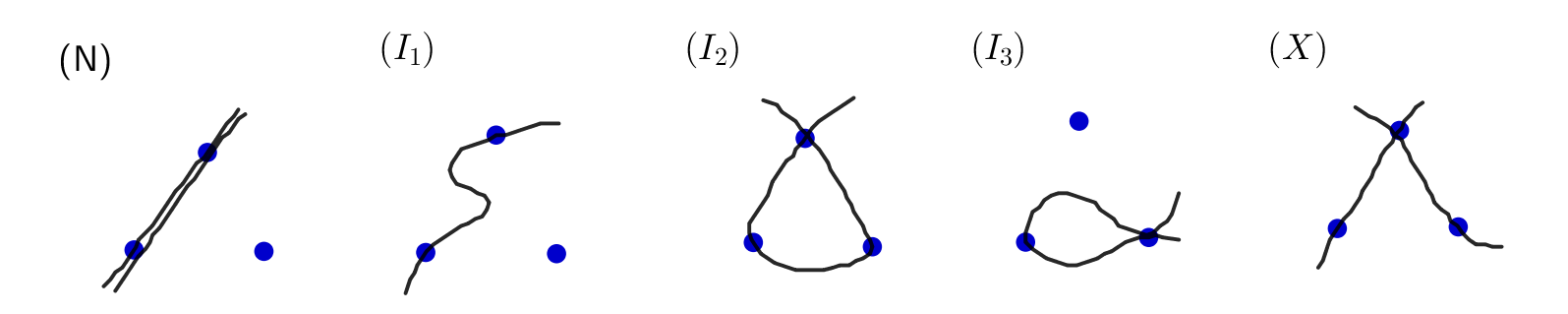}
\label{invcurves}
\end{figure}

\ni{\bf 6.1} Let $M$ be a fake projective plane with $K_M\equiv3L$ and $\Aut(M)=C_3$, where $L$ is invariant under $\Aut(M)$. Suppose that $\HH^0(M,2L)\neq0$ and $\Sigma\sim 2L$ is an $\Aut(M)$-invariant curve from Lemma \ref{invcurve1}. Let $(n,h^0(\delta),x)$ be the triple associated to $\Sigma_{\rm red}$ as in Lemma \ref{invirr}. According to Proposition \ref{FP} and Lemma \ref{invirr}, there is a finite list of possible $\Sigma$ according to the triple $(n,h^0(\delta),x)$:
\begin{enumerate}
\item[$(N)$] $(n,h^0(\delta),x)=(2,0,1)$: $\Sigma=2C$, where $C\equiv L$ is smooth and has two smooth fixed points;
	\item[$(I_1)$] $(n,h^0(\delta),x)=(2,0,2)$: $\Sigma$ is a smooth curve of $g(C)=6$ and has two smooth fixed points;
	\item[$(I_2)$] $(n,h^0(\delta),x)=(4,1,1)$: $\Sigma$ has one fixed node as the unique singularity and two smooth fixed points;
	\item[$(I_3)$] $(n,h^0(\delta),x)=(3,2,1)$: $\Sigma$ has one fixed tacnode as the unique singularity and one smooth fixed point. 	   
 \item[($X$)] $(n,h^0(\delta),x)=(4,1,2)$: $\Sigma=\Sigma_1+\Sigma_2$ has one fixed node $\{p\}=\Sigma_1\cap\Sigma_2$ as the unique singularity and one smooth fixed point along each $\Sigma_i$. Both $\Sigma_i$'s are $\Aut(M)$-invariant smooth curves of $g(\Sigma_i)=3$ with $x_i:=h^1(\Sigma_i,\cO_{\Sigma_i})^{\rm inv}=1$.
\end{enumerate} 
Only the case $(X)$ needs to be explained: If $\Sigma=\Sigma_1+\Sigma_2$, then $\Sigma_i$'s are smooth of genus 3 by Lemma \ref{Sch} and invariant under $C_3$. Apply Lemma \ref{LFP}, we see there are two fixed points on each $\Sigma_i$ and hence the description.\footnote{Each $\Sigma_i$ has $x_i=1$ and hence $E:=\Sigma_i/C_3\subseteq M/C_3$ is an elliptic curve. If $E\hookrightarrow M/C_3$ lifts to $E\rightarrow M$, then this contradicts the hyperbolicity of $M$. However, as $M/C_3$ is singular, the lifting does not always exist.} We will show in Lemma \ref{X} that case $(X)$ does not happen and hence the picture of all invariant curves is as in Figure \ref{invcurves}.

\ms

\ni{\bf 6.2} We will study $\Aut(M)$-invariant curves by investigating their geometry on $Y:=M/\Aut(M)$ and its minimal resolution. Again as in Section \ref{Sec4} and \ref{Sec5}, we want to show the absence of these invariant curves. The following structure theorem on $Y$ is crucial.
\begin{theo}[\cite{CS,K1}] Suppose that a fake projective plane $M$ has $\Aut(M)=C_3$. Then the surface $Y=M/\Aut(M)$ has 3 singularities of type $\frac{1}{3}(1,2).$
\end{theo}
The quotient surface $Y=M/\Aut(M)$ is a $\bQ$-homology projective plane and say $\Sing(Y)=\{a,b,c\}$. By abuse of notation, we also denote $\{a,b,c\}\subseteq M$ the $\Aut(M)$-fixed points corresponding to ${\rm Sing}(Y)$. Let  $\mu:Z\rightarrow Y$ be the minimal resolution. Then $Z$ is a minimal surface of general type with $p_g(Z)=q(Z)=0$ and $K_Z^2=3.$ Let $C\equiv kL$ be an $\Aut(M)$-invariant curve. We study the quotient curves $C_Y=C/C_3$ on $Y$ and its proper transform $C_Z:=\mu^{-1}_*(C_Y)$ on the minimal resolution $Z$of $Y$. Note that $$g(C_Z^\nu)=\HH^0(C^\nu,\omega_{C^\nu})^{\rm inv}=x.$$

\begin{lemm}\label{7} The Picard number of $Z$ is $\rho(Z)=7$. 
\end{lemm}
\ni{\bf Proof.} Noether equality implies that $\chi_{\rm top}(Z)=9.$ By Hodge decomposition and $p_g(Z)=0$, we have $h^{1,1}(Z)=b_2(Z)=7$. Since $p_g(Z)=q(Z)=0$, by exponential sequence $\rho(Z)=h^{1,1}(Z)=7$. 
\qed

\ms

\begin{lemm}\label{intnum} Over each singular point $\star\in{\rm Sing}(Y)$, $\mu^{-1}(\star)=E_{\star1}\cup E_{\star2}$ is a simple normal crossing divisor with $E_{\star 1}^2=E_{\star2}^2=-2$ and $E_{\star 1}\cdot E_{\star2}=1$. Write $\mu^*C_Y=C_Z+E_Z$, where $E_Z$ is $\mu$-exceptional. Locally at $\star\in{\rm Sing}(Y)$, we have 
$$(E_Z,E_Z^2)=\begin{cases}(E_{\star 1}+E_{\star 2},-2)\ &{\rm if}\ \star\in C_Y\ {\rm is\ nodal},\\ 
(\frac{2}{3}E_{\star 1}+\frac{1}{3}E_{\star 2},-\frac{2}{3})\ &{\rm if}\ \star\in C_Y\ {\rm is\ smooth\ and}\ C_Z\cap E_{\star1}\neq\emptyset,\\
(\frac{4}{3}E_{\star 1}+\frac{2}{3}E_{\star 2},-\frac{8}{3})\ &{\rm if}\ \star\in C_Y\ {\rm is\ a\ tacnode\ and}\ C_Z\cap E_{\star1}\neq\emptyset.\\
\end{cases}$$
\end{lemm}
\ni{\bf Proof.} A $\frac{1}{3}(1,2)$-point is an $A_2$ rational double point and can be described as the point \{$O=(0,0,0)\in S:=(XY-Z^3=0)\subseteq \bC^3$ via 
$$\bC^2\twoheadrightarrow\bC^3/C_3\xrightarrow{\cong}S, \ (x,y)\mapsto(X,Y,Z)=(x^3,y^3,xy).$$
Moreover, the proper transform $\tilde{S}$ of $S$ in $\Bl_O\bC^3$ is the minimal resolution of $O\in S$. The exceptional divisor of $\mu:\tilde{S}\rightarrow S$ in $\bC^3\times P^2_\bC$ with $[u:v:w]$ the homogeneous coordinate is given by $P_\bC^1\cup P_\bC^1\cong(uv=0)\subseteq \{O\}\times P_\bC^2\subseteq\Bl_O\bC^3.$ The rest is easy and left to the reader.
\qed

\ms

As $K_Z=\mu^*K_Y$, $K_Z^2=3$, and $C_Y^2=\frac{k^2}{3}$, we get
\begin{align*}
p_a(C_Z)=1+\frac{1}{2}(K_Y.C_Y+C_Y^2-C_Z.E_Z)
=1+\frac{k(k+3)}{6}-\frac{1}{2}C_Z\cdot E_Z.
\end{align*} 
We can now list the relevant intersection numbers on $Z$, see Table \ref{intnumZ}. Here for $C_Z$ we take $C=\Sigma_{\rm red}$ in the case $(N)$ and $C=\Sigma$ for all the other cases.
\begin{table}
\caption{Invariant curves on $Z$}
\begin{center}
	\begin{tabular}{|c|c|c|c|c|c|c|}\hline
		& $(N)$ & $(I_1)$ & $(I_2)$ & $(I_3)$ & $(X)$   \\ \hline
		$C_Z\cdot E_Z=-E_Z^2$ & $\frac{4}{3}$ & $\frac{4}{3}$ & $\frac{10}{3}$ & $\frac{10}{3}$ & $\frac{10}{3}$  \\ \hline
		$C_Z^2=\frac{k^2}{3}-C_Z\cdot E_Z$ & -1 & 0 & -2 & -2 & -2  \\ \hline
		$p_a(C_Z)$ & 1 & 2 & 1 & 1 & 1 \\ \hline
		$K_Z.C_Z$ & 1 & 2 & 2 & 2 & 2   \\ \hline
		$g(C^\nu_Z)=x$ & 1 & 2 & 1 & 1 & 2  \\ \hline
	\end{tabular}
\end{center}	
	\label{intnumZ}
\end{table}

We first show that the case $(X)$ does not happen after the following observation.
\begin{lemm}\label{invtor2} Let $M$ be a fake projective plane with $\Aut(M)=C_3$. Then $12\tau=0$ for any $\Aut(M)$-invariant torsion line bundle $\tau$.
\end{lemm}
\ni{\bf Proof.} This follows from the proof of Lemma \ref{invtor} that $\Aut(M)$-invariant line bundles always descend to $M/\Aut(M)$. From the list of $\HH_1(M/\Aut(M),\bZ)$ in the file \href{http://www.maths.usyd.edu.au/u/donaldc/fakeprojectiveplanes/registerofgps.txt}{registerofgps.txt} from the weblink of \cite{CS}, we see that $12\tau=0$, cf. Table \ref{Tab:Main} and \ref{Tab:unproven}.
\qed

\begin{rem} For any fake projective plane with $\Aut(M)=C_3$ in Table \ref{Tab:Main}, $6\tau=0$ for any $\Aut(M)$-invariant torsion line bundle $\tau$. 
However, we still use Lemma \ref{invtor2} for the discussion on invariant curves: the following Lemma \ref{X} illustrates that our approach potentially can 
work for any fake projective plane with non-trivial automorphisms. Hence in Section \ref{Sec6} to \ref{Sec8}, all the statement are given in its most general form intensionally for future reference. 
\end{rem}

\begin{lemm}\label{X} Let $M$ be a fake projective plane with $\Aut(M)=C_3$ and  $L$ is an $\Aut(M)$-invariant numerical cubic root of $K_M$, i.e., $K_M\equiv 3L$. Then an invariant curve $\Sigma\sim 2L$ of type $(X)$ does not exist.
\end{lemm}
\ni{\bf proof:} Let $\Sigma=\Sigma_1+\Sigma_2\sim 2L$. Consider the image curves $\Sigma_i^Y:=\pi(\Sigma_i)$ on $Y=M/\Aut(M)$, where $\pi:M\rightarrow Y$ is the quotient map, and their proper transforms $\Sigma_i^Z:=\mu^{-1}(\Sigma_i^Y)$, $i=1,2$, on the minimal resolution $Z$. Assume that $\{a,b\}\subseteq\Sigma_1$ and $\{a,c\}\subset\Sigma_2$, where we have identified ${\rm Fix}(M)=\Sing(Y)$.  By Lemma \ref{intnum}, we may assume that  
$$\begin{cases} \mu^*\Sigma_1^Y=\Sigma_1^Z+\frac{2}{3}E_{a1}+\frac{1}{3}E_{a_2}+\frac{2}{3}E_{b1}+\frac{1}{3}E_{b_2}\\
\mu^*\Sigma_2^Y=\Sigma_2^Z+\frac{1}{3}E_{a1}+\frac{2}{3}E_{a_2}+\frac{2}{3}E_{c1}+\frac{1}{3}E_{c_2}
\end{cases}.$$
Since $12\tau=0$ by Lemma \ref{invtor2} and $\Sigma_i$'s are $\Aut(M)$-invariant, associated to the three linearly independent\footnote{This can be checked easily by considering the vanishing order along $\Sigma_1$ and $\Sigma_2.$} divisors $12\cdot(2\Sigma_1)\sim12\cdot(2\Sigma_2)\sim12\cdot(\Sigma_1+\Sigma_2)$ in $|24L|$, we find the following three linearly independent elements in $|24\mu^*L_Y-8(E_{a1}+E_{a_2})|$:
$$\begin{cases} S:=24\Sigma_1^Z+8E_{a1}+16E_{b1}+8E_{b_2},\\
T:=24\Sigma_2^Z+8E_{a_2}+16E_{c1}+8E_{c_2},\ {\rm and}\\
U:=12(\Sigma_1^Z+\Sigma_2^Z)+4E_{a1}+4E_{a_2}+8E_{b1}+4E_{b_2}+8E_{c1}+4E_{c_2}.
\end{cases}$$
In particular, the subsystem $\Lambda:=\left<S,T,U\right>$ has a unique base point at $z:=E_{a1}\cap E_{a2}$ (of length $S\cdot T=64$). Let $\mu':Z'\rightarrow Z$ be the blow up at $z\in Z$ with the unique exceptional divisor $E$. Then $\Lambda'=\left<S',T',U'\right>\subseteq|\mu'^*(24\mu^*L_Y-8(E_{a1}+E_{a_2}))-8E|$ with 
$$\begin{cases} S':=24\Sigma_1^{Z'}+8E_{a1}+16E_{b1}+8E_{b_2},\\
T':=24\Sigma_2^{Z'}+8E_{a_2}+16E_{c1}+8E_{c_2},\ {\rm and}\\
U':=12(\Sigma_1^{Z'}+\Sigma_2^{Z'})+4E_{a1}+4E_{a_2}+8E_{b1}+4E_{b_2}+8E_{c1}+4E_{c_2},
\end{cases}$$
is base point free. Here we abuse the notion by denoting $E_{\star i}$ again its proper transform on $Z'$. There is a morphism $\varphi:=\varphi_{\Lambda}:Z'\rightarrow \mathbb{P}(\Lambda)\cong P_\bC^2$, where $(S')^2=0$ implies that $\dim(\varphi(Z))=1$. As the image is dominated by $E\cong P^1_\bC$, the Stein factorization of $\varphi$ induces a morphism $\varphi':Z\rightarrow P_\bC^1$ with connected fibers. Since $S'$ and $T'$ are connected, they are two special fibers of $\varphi'$. This clearly is impossible since $U'$ shares some common components with $S'$ and $T'$.
\qed

\ms

We remark that to rule out pairs of invariant curves, the construction of a free pencil as in the proof of Lemma \ref{X} will appear several times in different occasions.

\ms

\ni{\bf 6.3} We are not able to directly show that $\HH^0(M,2L)=0$ as required in Lemma \ref{ECeasy} by ruling out all types of invariant curves from the list in {\bf 6.1} as done in Theorem \ref{1} and \ref{2}. Instead, we investigate how different $\Aut(M)$-invariant curves intersect if there are many numerical cubic roots of $K_M$ with non-vanishing cohomology. Suppose that a line bundle $L'\ncong L$ is another $\Aut(M)$-invariant generator of $\NS(M)_\bQ$. We assume that both $\HH^0(M,2L)\neq0$ and $\HH^0(M,2L')\neq0$. Let $\Sigma\sim2L$ and $\Sigma'\sim2L'$ be two $\Aut(M)$-invariant curves from Lemma \ref{invcurve1}. The following proposition shows that the intersection behavior of $\Sigma$ and $\Sigma'$ is rather restricted. 

\begin{prop}\label{double} In the setting above, the two curves $\Sigma$ and $\Sigma'$ can only intersect along $\Aut(M)$-fixed points of $M$.
\end{prop}

\ni{\bf Proof.} We will prove that the set of divisors  $\{\Sigma_Z,\Sigma'_Z,E_{a1}, E_{a2},E_{b1}, E_{b2},E_{c1}, E_{c2}\}$ is linearly independent in $N^1(Z)$ if $\Sigma$ and $\Sigma'$ intersect at a non-$\Aut(M)$-fixed point. This will contradict Lemma \ref{7}. Clearly, it is enough to show that the 8 by 8 intersection matrix $I$ of this set of eight curves has non-zero determinant.

Apriori there are 16 possible types of $(\Sigma,\Sigma')$ from the list of $(N)$, $(I_1)$, $(I_2)$, and $(I_3)$. Each pair we have to consider possible intersection figuration to get $\Sigma\cdot\Sigma'=4.$ Recall that $\{a,b,c\}\subseteq M$ is the set of $\Aut(M)$-fixed points corresponding to ${\rm Sing}(Y)$.

Suppose now $\Sigma$ and $\Sigma'$ intersect at a point $o'\notin\{a,b,c\}$. Hence they intersect along the $C_3$-orbit $C_3\cdot o'$. But then $\Sigma$ and $\Sigma'$ must intersect transversally at $o'$ by 
$$4=\Sigma\cdot\Sigma'\geq3{\rm mult}_{o'}(\Sigma\cap\Sigma').$$
In particular, $\Sigma$ and $\Sigma'$ intersect transversally at one another $\Aut(M)$-fixed point of $M$, say at $a$. Note that this does not happen if one of $\Sigma$ and $\Sigma'$ is of type $(N)$: otherwise 
$$\Sigma\cdot\Sigma'\geq2\Sigma_{\rm red}\cdot\Sigma'_{\rm red}\geq 2\cdot 3.$$ 
Also, none of $\Sigma$ and $\Sigma'$ is of type $(I_2)$, otherwise they must intersect along at least two $\Aut(M)$-fixed points. 

Hence up to reordering, the intersection matrix $I$ can only be one of the following possibilities, where we have used Lemma \ref{intnum} and Table \ref{intnumZ}. Note that by assumption, $\Sigma_Z\cdot\Sigma'_Z=3$.

{\bf Case 1}: $(\Sigma,\Sigma')$ is of type $(I_1,I_1)$. Say $\{a,b\}\subseteq\Sigma$ and $\{a,c\}\subseteq\Sigma'$ with $\Sigma_Z$ intersecting $E_{a1}$ and $E_{b1}$ while $\Sigma'_Z$ intersecting $E_{a2}$ and $E_{c1}$:
\begin{align*} I=\left(\begin{array}{cccccccc}
0 & 3 & 1 & 0 & 1 & 0 & 0 & 0 \\
3 & 0& 0 & 1 & 0 & 0 & 1 & 0 \\
1 & 0 & -2& 1 & 0 & 0 & 0 & 0 \\
0 & 1 & 1 &-2 & 0 & 0 & 0 & 0 \\
1 & 0 & 0 & 0 &-2 & 1 & 0 & 0 \\
0 & 0 & 0 & 0 & 1 &-2 & 0 & 0 \\ 
0 & 1 & 0 & 0 & 0 & 0 &-2 & 1 \\
0 & 0 & 0 & 0 & 0 & 0 & 1 & -2
\end{array}\right)\ {\rm and}\ \det(I)=-252.
\end{align*}
{\bf Case 2}: $(\Sigma,\Sigma')$ is of type $(I_1,I_3)$. Say $\{a,b\}\subseteq\Sigma$ and $\{a,c\}\subseteq\Sigma'$. Suppose that $\Sigma_Z$ intersects $E_{a1}$ and $E_{b1}$ and $\Sigma'_Z$ intersects $E_{a2}$ and $E_{c1}$:
\begin{align*} I=\left(\begin{array}{cccccccc}
0 & 3 & 1 & 0 & 1 & 0 & 0 & 0 \\
3 & -2& 0 & 1 & 0 & 0 &2 & 0 \\
1 & 0 & -2& 1 & 0 & 0 & 0 & 0 \\
0 & 1 & 1 &-2 & 0 & 0 & 0 & 0 \\
1 & 0 & 0 & 0 &-2 & 1 & 0 & 0 \\
0 & 0 & 0 & 0 & 1 &-2 & 0 & 0 \\ 
0 & 2 & 0 & 0 & 0 & 0 &-2 & 1 \\
0 & 0 & 0 & 0 & 0 & 0 & 1 & -2
\end{array}\right)\ {\rm and}\ \det(I)=-252.
\end{align*}
{\bf Case 3}:  $(\Sigma,\Sigma')$ is of type $(I_3,I_3)$. Say $\{a,b\}\subseteq\Sigma$ and $\{a,c\}\subseteq\Sigma'$. Suppose that $\Sigma_Z$ intersects $E_{a1}$ and $E_{b1}$ and $\Sigma'_Z$ intersects $E_{a2}$ and $E_{c1}$:
\begin{align*} I=\left(\begin{array}{cccccccc}
-2 & 3 & 1 & 0 & 2 & 0 & 0 & 0 \\
3 & -2& 0 & 1 & 0 & 0 & 2 & 0 \\
1 & 0 & -2& 1 & 0 & 0 & 0 & 0 \\
0 & 1 & 1 &-2 & 0 & 0 & 0 & 0 \\
2 & 0 & 0 & 0 &-2 & 1 & 0 & 0 \\
0 & 0 & 0 & 0 & 1 &-2 & 0 & 0 \\ 
0 & 2 & 0 & 0 & 0 & 0 &-2 & 1 \\
0 & 0 & 0 & 0 & 0 & 0 & 1 & -2
\end{array}\right)\ {\rm and}\ \det(I)=-252.
\end{align*}
Since all the determinants are non-zero, this proves the result as discussed before. 
\qed

\ms

\ni{\bf 6.4} In this subsection, we provide the complete list of all possible configurations of two invariant curves $\Sigma\neq\Sigma'$ in the numerical class $2L$. We first classify the local intersection configurations. 
\begin{lemm}\label{localint}Suppose that there are two $\Aut(M)$-invariant curves $\Sigma$ and $\Sigma'$ of numerical type $2L$ intersecting at a fixed point $a\in M$. Then $\Sigma$ and $\Sigma'$ share no common component. Moreover, up to relabelling,  local analytically around the point $a$ it is in one of the following configurations:
			\begin{enumerate}
				\item If both $\Sigma_{\rm red}$ and $\Sigma'_{\rm red}$ are uni-branched, then $\Sigma_{\rm red}\cup\Sigma'_{\rm red}$ is in one of the following forms:
$$
\begin{array}{|c|c|c|}
\hline
{\rm Notation}&{\rm local\ equation} &{\rm mult}_a(\Sigma_{\rm red}\cap \Sigma'_{\rm red})\\ \hline\hline
(tr) & xy=0 &1 \\ \hline
(tan-sm) & x(x-y^2)=0 &2\\ \hline
(tan-tan) & (x-y^2)(x+y^2)=0 &4 \\ \hline
\end{array}
$$									
		\item If $\Sigma_{\rm red}$ is uni-branched but $\Sigma'=\Sigma'_{\rm red}$ is two-branched at $a$, then $\Sigma_{\rm red}\cup\Sigma'_{\rm red}$ is in one of the following forms: 
$$
\begin{array}{|c|c|c|}
\hline
{\rm Notation}&{\rm local\ equation} &{\rm mult}_a(\Sigma_{\rm red}\cap \Sigma'_{\rm red})\\ \hline\hline
(tr-tac) & y\cdot x(x-y^2)=0 &2\\ \hline
(tan-node) & (x-y^2)\cdot xy=0 &3 \\ \hline
(tan-tac) & x\cdot(x^2-y^4)=0 &4 \\ \hline
\end{array}
$$					
	\end{enumerate}
\begin{figure}
\centering
\caption{LOCAL INTERSECTIONS}
\includegraphics*[scale=0.5]{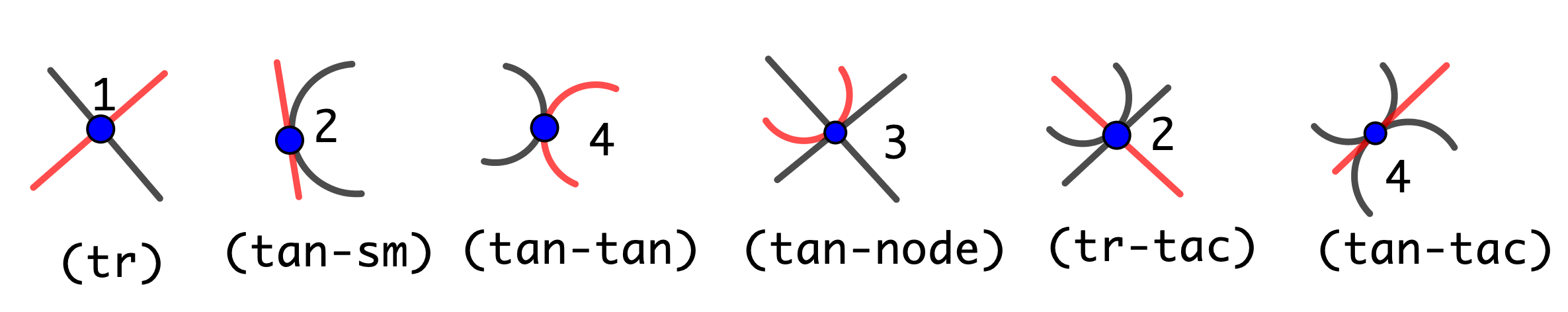}
\end{figure}

\end{lemm} 
\ni{\bf Proof.} We follow the computation in Lemma \ref{local}. Note that $\Sigma\cdot\Sigma'=(2L)^2=4$ and the intersection multiplicity satisfies $\mult_p(\Sigma\cap\Sigma')\leq4,\ \forall\ p\in\Sigma\cap\Sigma',$
unless they share a common component. 

Suppose now $\Sigma$ and $\Sigma'$ share a common branch near $a$. Then as $\Sigma$ and $\Sigma'$ are algebraic curves, they share an irreducible component. But then $\Sigma$ or $\Sigma'$ has to be reducible as $\Sigma\neq\Sigma'.$ Hence $\Sigma$ or $\Sigma'$ has to be of type $(X)$, which violates Lemma \ref{X}. Hereafter we assume that $\Sigma$ and $\Sigma'$ share no common branch near $a$. 

First observe that locally at $a\in M$, one of $\Sigma_{\rm red}$ and $\Sigma'_{\rm red}$ must be uni-branched: Suppose that both of them are two-branched. In particular, $\Sigma$ and $\Sigma'$ are reduced from Lemma \ref{FP}. Assume that $\Sigma=(xy=0)$ in an analytic neighborhood of $a$ of weight $(x,y)=(1,2)$. If $\Sigma'$ is also nodal at $a$, then from Lemma \ref{local} a branch of $\Sigma'$ is of the form $(x\pm\epsilon y^2+{\rm h.o.t}=0)$ or $(y\pm \epsilon x^2+{\rm h.o.t}=0)$, where $\epsilon\in\{0,1\}$. If $\Sigma'=(xy=0)$, then $\Sigma=\Sigma'$ local analytically and $\Sigma=\Sigma'$ on $M$ (for being algebraic curves with non-isolated intersection). If $\Sigma'\neq(xy=0)$, then the intersection multiplicity $\mult_a(\Sigma\cap\Sigma')>4$, which is absurd. 
Suppose now that $\Sigma=(x(x-y^2)=0)$ is a tacnode at $a$. By the same consideration as above, either we get $\Sigma'=\Sigma$ or $\mult_a(\Sigma\cap\Sigma')>4$, which is again impossible. 

We now classify possible local intersection configuration in two cases.

\ni{\bf Case 1}: Both $\Sigma_{\rm red}$ and $\Sigma'_{\rm red}$ are uni-branched at $a$. 

Since $\Sigma_{\rm red}\neq\Sigma'_{\rm red}$ near $a$ , by Lemma \ref{local} it is easy to see that  $\Sigma_{\rm red}\cup\Sigma'_{\rm red}$ local analytically is in one of the following form :
\begin{enumerate}
	\item $((xy=0)$ with  ${\rm mult}_a(\Sigma_{\rm red}\cap \Sigma'_{\rm red})=1;$
	\item $(x(x-y^2)=0)$ with  ${\rm mult}_a(\Sigma_{\rm red}\cap \Sigma'_{\rm red})=2;$ 
	\item $((x-y^2)(x+y^2)=0)$ with  ${\rm mult}_a(\Sigma_{\rm red}\cap \Sigma'_{\rm red})=4.$ 
\end{enumerate}
The first case is when $\Sigma_{\rm red}$ and $\Sigma'_{\rm red}$ intersect transversally $a$, while the last two cases are when  $\Sigma_{\rm red}$ and $\Sigma'_{\rm red}$ intersection tangentially at $a$.

\ni{\bf Case 2}: $\Sigma$ is uni-branched at $a$ but $\Sigma'=\Sigma'_{\rm red}$ is two-branched at $a$.

From the list in {\bf 6.1},  $\Sigma'=\Sigma'_{\rm red}$ if it is two-branched at a fixed point. There are two cases.

{\bf Subcase 2.1:} $a\in\Sigma'$ is a node.

Since there are only two eigen-directions at a fixed point $a\in M$, $\Sigma$ must intersect $\Sigma'$ tangentially at $a$. Say locally $\Sigma'=(xy=0)$ with weight of $(x,y)$ being $(1,2)$.  From $\mult_a(\Sigma\cap\Sigma')\leq4$, $\Sigma_{\rm red}$ can only be $(x-y^2+{\rm h.o.t.}=0)$ or $(y-x^2+{\rm h.o.t.}=0)$. 
Hence near $a$ we have $\Sigma_{\rm red}\cup\Sigma'\cong (xy(x-y^2)=0)$ with $\mult_a(\Sigma_{\rm red}\cap\Sigma')=3$. But then $\Sigma=\Sigma_{\rm red}$.

{\bf Subcase 2.2.a:} $a\in\Sigma'$ is a tacnode and $\Sigma$ intersects $\Sigma'$ transversally at $a$.

We may assume that $\Sigma'$ locally near $a$ has the equation $x(x-y^2)=0$. Since the intersection is transversal, 
we have $\Sigma_{\rm red}\cup\Sigma'\cong (yx(x-y^2)=0)$ and $\mult_a(\Sigma_{\rm red}\cap\Sigma')=2$. Note that it is possible $\Sigma$ to be of type $(N)$ in this case.

{\bf Subcase 2.2.b:} $a\in\Sigma'$ is a tacnode and $\Sigma$ intersects $\Sigma'$ tangentially at $a$.

Assume that $\Sigma'$ locally near $a$ has the equation $x(x-y^2)=0$. Since $\Sigma$ and $\Sigma'$ share no common component, local equation of $\Sigma_{\rm red}$ near $a$ is of the form 
$x\pm\epsilon y^2+{\rm h.o.t.}$ with $\epsilon\in\{0,1\}$. In particular, $\mult_a(\Sigma\cap\Sigma')\geq4$, and equality holds only if locally we have $\Sigma_{\rm red}\cup\Sigma'\cong(x(x^2-y^4)=0)$.
\qed

\ms
Now we classify possible intersection configurations of two different invariant curves. In the following, the intersection configurations refer to the terminology given in Lemma \ref{localint}.  An expression such as $3(tan-node)+1(tr)$ reflects that the intersection multiplicity of the two curves is $4$ with $3$ contributed by an intersection configuration $tan-node$ and $1$ contributed by an intersection configuration $tr$ given by the reduced parts of the curves at the corresponding intersection 
points.

\begin{lemm}\label{2Inv} Given two distinct $\Aut(M)$-invariant curves $\Sigma$ and $\Sigma'$ in the numerical class $2L$. Then $|\Sigma\cap\Sigma'|\leq2$ and they share no common component. Up to relabelling, the type of $(\Sigma,\Sigma')$ is in one of the following pairs:
	\begin{enumerate}		
		\item[(1a)] $(I_1,I_2)$ with intersection configuration $3(tan-node)+1(tr)$;
		\item[(1b-1)] $(N,I_3)$ with intersection configuration $2\cdot2(tr-tac)$;
		\item[(1b-2)] $(I_3,I_3)$ with intersection configuration $2(tr-tac)+2(tr-tac)$;
		\item[(1b-3)] $(I_1,I_3)$, $(I_2,I_3)$ with intersection configuration $2(tr-tac)+2(tan-sm)$;
		\item[(1c)] $(I_1,I_3)$, $(I_3,I_3)$ with intersection configuration $4(tan-tac)$;
		\item[(2a)] $(I_1,I_1)$, $(I_1,I_3)$, $(I_3,I_3)$ with intersection configuration $4(tan-tan)$;
		\item[(2b)] $(I_1,I_1)$, $(I_1,I_2)$ of  intersection type $2(tan-sm)+2(tan-sm)$;
		\item[(3a)] $(N,N)$ with intersection configuration $4\cdot 1(tr)$;
		\item[(3b)] $(N,I_1)$, $(N,I_2)$ with intersection configuration $2\cdot(1(tr)+1(tr))$;
		\item[(3c)] $(N,I_1)$, $(N,I_3)$ with intersection configuration $2\cdot2(tan-sm)$.
	\end{enumerate}
	     
\end{lemm}
 
\ni{\bf Proof.} By Proposition \ref{double}, $\Sigma$ and $\Sigma'$ only intersect along $\Aut(M)$-fixed points. By Lemma \ref{localint}, $\Sigma$ and $\Sigma'$ share no common component. If $\Sigma\cap\Sigma'=\{a,b,c\}$, then both $\Sigma$ and $\Sigma'$ possess nodes. It is clear in this case $\Sigma\cdot\Sigma'>4$ from Lemma \ref{localint}, which is absurd. 

Hereafter we assume that $\Sigma$ and $\Sigma'$ share no common component and $a\in\Sigma\cap\Sigma'$. By Lemma \ref{localint}, we study their intersection configuration by considering the following two cases (possibly after relabelling):
\begin{enumerate}
	\item $\Sigma'=\Sigma'_{\rm red}$ is two-branched at $a$.
	\item $\Sigma_{\rm red}$ and $\Sigma'_{\rm red}$ are uni-branched (and hence smooth) at all intersection points.
\end{enumerate}
For simplicity, denote by $m_p=\mult_p(\Sigma\cap\Sigma')$ for $p=a,b$, or $c$. 

{\bf Case (1a):} $a\in\Sigma'$ is a node.

From the list of invariant curves, $\Sigma'$ is of type $(I_2)$. By Lemma \ref{localint}, $\Sigma$ intersects tangentially at $a\in\Sigma'$ with $m_a=3$ and hence cannot be of type $(N)$. Moreover, $\Sigma$ and $\Sigma'$ can only intersect transversally at another fixed point, say $b\in{\rm Fix}(M)$. In particular, $\Sigma$ has two smooth fixed points and is one of $(I_1)$ or $(I_2)$. The last case is impossible since then these two curves intersect at all three fixed points and $\Sigma\cdot\Sigma'>4.$ 

{\bf Case (1b):} $a\in\Sigma'$ is a tacnode.

Hence $\Sigma'$ is of type $(I_3)$ and $m_a=2$ or $4$ by Lemma \ref{localint}. 

Suppose that $m_a=2$. Note that $a\in\Sigma_{\rm red}$ is smooth. If $\Sigma$ is of type $(N)$, then $\Sigma$ and $\Sigma'$ intersect transversally at the unique intersection point $a$. This is type $(N, I_3)$ in  (1b-1). We may assume now $\Sigma$ is reduced. Since $\Sigma'$ has only one more smooth fixed point $b$, we must have $m_b=2$. If $b\in\Sigma$ is a tacnode, then they intersect transversally at $b$ and we have type $(I_3,I_3)$ as (1b-2). If $b\in\Sigma$ is smooth, then they intersect tangentially at $b$ and $\Sigma$ has at least two smooth fixed points. Hence $(\Sigma,\Sigma')$ can be of type $(I_1,I_3)$ or $(I_2,I_3)$ as in (1b-3). 

If $m_a=4$, then $a\in\Sigma=\Sigma_{\rm red}$ is smooth and is the unique intersection point. Hence $\Sigma$ can only be of type $(I_1)$ or $(I_3).$ This is case $(1c)$.

We assume now that $\Sigma$ and $\Sigma'$ share no common components, reduced, and are uni-branched (and hence smooth) at all fixed points. In particular, none of them is of type $(N)$ and the local configuration is as in Lemma \ref{localint} (1).

{\bf Case (2a):} $\mult_p(\Sigma\cap\Sigma')= 4$ for some $p\in\Sigma\cap\Sigma'$.

Since $\Sigma\cdot\Sigma'=4$, we may assume that $a\in\Sigma\cap\Sigma'$ is the unique intersection point. As each invariant curve has at least one fixed point, $\Sigma$ and $\Sigma'$ both must have exactly two fixed points. Since none of them are of type $(N)$, the only possible types of $(\Sigma,\Sigma')$ are $(I_1,I_1)$, $(I_1,I_3)$, and $(I_3,I_3).$

{\bf Case (2b):} $\mult_p(\Sigma\cap\Sigma')= 2$ at two fixed points $p\in\{a,b\}$. 

Clearly, $\Sigma\cap\Sigma'=\{a,b\}$. As none of two curves are of type $(N)$, $\Sigma$ and $\Sigma'$ are smooth along $a,b$ and intersect tangentially at both $a$ and $b$. In particular, none of them is of type $(I_3).$ Since they cannot simultaneously have three fixed points, the remaining possible types are $(I_1,I_1)$ and $(I_1,I_2)$.

{\bf Case (2c):} $\mult_a(\Sigma\cap\Sigma')= 2$ at exactly one fixed point.

Then $\Sigma$ and $\Sigma'$ must intersect transversally at the other two fixed points. In particular, both of them have three fixed points and can only be of type $(I_2)$. This violates the assumption that both of them must be smooth at all fixed points. 

{\bf Case (2d):} $\mult_p(\Sigma\cap\Sigma')= 1$ for all $p\in\Sigma\cap\Sigma'$.

Since each invariant curve can have at most two smooth fixed points, it is impossible in this case to have $\Sigma\cdot\Sigma'=4.$

For the remaining cases, we assume that one of them is of type $(N)$.

{\bf Case (3a):} Both $\Sigma=2C$ and $\Sigma'=2C'$ are of type $(N)$.

Since $C\cdot C'=\frac{1}{4}\Sigma\cdot\Sigma'=1$, $C$ and $C'$ intersect transversally at a fixed point. 

{\bf Case (3b):} $\Sigma=2C$ and $\mult_p(C\cap\Sigma')= 1$ for all $p\in C\cap\Sigma'$.

Then $\Sigma'$ must have two smooth fixed points and this can happen if it is of type $(I_1)$ or $(I_2)$.

{\bf Case (3c):} $\Sigma=2C$ and $\mult_a(C\cap\Sigma')=2$ at a fixed point $a$.

Then $\Sigma\cap\Sigma'=\{a\}$. By Lemma \ref{localint}, $C$ and $\Sigma'$ intersect tangentially at a smooth point $a\in\Sigma'$ or $C$ and $\Sigma'$ intersect transversally at a tacnode $a\in \Sigma'$. If we are in the former case, then $\Sigma'$ cannot have three fixed points and hence is of type $(I_1)$ or $(I_3)$. In the latter case, we can only have type $(N,I_3)$, which duplicates with case (1b-1).
\qed

\ms

\ni{\bf 6.5} By applying a similar argument for proving Proposition \ref{double} and Lemma \ref{X}, we are able to rule out some cases in Lemma \ref{2Inv}. 

\begin{lemm}\label{m=4} The case (2a) (of two curves intersecting at exactly one point of multiplicity four) in Lemma \ref{2Inv} does not occur.
\end{lemm}
\ni{\bf Proof:} We follow the same proof of Lemma \ref{X} by showing that the intersection matrix $I$ of the eight curves $\{\Sigma_Z,\Sigma'_Z,E_{a1}, E_{a2},E_{b1}, E_{b2},E_{c1}, E_{c2}\}$ is non-degenerate. It follows that $\rho(Z)\geq8$, which violates Lemma \ref{7}. 

In (2a), we have $(I_1,I_1)$, $(I_1,I_3)$ or $(I_3,I_3)$, and say with $m_a=4$. From the local description in Lemma \ref{intnum} and Lemma \ref{localint}, it is easy to see that $\Sigma_Z\cdot\Sigma_Z'=0$. We assume on $Z$ that $\Sigma_Z$ intersects $E_{a1}$ and $E_{c1}$, while $\Sigma'_Z$ intersects $E_{a1}$ and $E_{b1}$. From Table \ref{intnumZ}, the intersection matrices respectively are 
\begin{align*} I(I_1,I_1)=\left(\begin{array}{cccccccc}
0 & 0 & 1 & 0 & 0& 0 &1 & 0 \\
0 &  0& 1 & 0 & 1 & 0 & 0 & 0 \\
1 & 1 & -2& 1 & 0 & 0 & 0 & 0 \\
0 & 0 & 1 &-2 & 0 & 0 & 0 & 0 \\
0 & 1 & 0 & 0 &-2 & 1 & 0 & 0 \\
0 & 0 & 0 & 0 & 1 &-2 & 0 & 0 \\ 
1 & 0 & 0 & 0 & 0 & 0 &-2 & 1 \\
0 & 0 & 0 & 0 & 0 & 0 & 1 & -2
\end{array}\right)\ {\rm with}\ \det(I)=36;
\end{align*}
\begin{align*}
I(I_1,I_3)=\left(\begin{array}{cccccccc}
0 & 0 & 1 & 0 & 0 & 0 & 1 & 0 \\
0 &  -2& 1 & 0 & 2 & 0 & 0 & 0 \\
1 & 1 & -2& 1 & 0 & 0 & 0 & 0 \\
0 & 0 & 1 &-2 & 0 & 0 & 0 & 0 \\
0 & 2 & 0 & 0 &-2 & 1 & 0 & 0 \\
0 & 0 & 0 & 0 & 1 &-2 & 0 & 0 \\ 
1 & 0 & 0 & 0 & 0 & 0 &-2 & 1 \\
0 & 0 & 0 & 0 & 0 & 0 & 1 & -2
\end{array}\right)\ {\rm with}\ \det(I)=36;
\end{align*}
\begin{align*}
I(I_3,I_3)=\left(\begin{array}{cccccccc}
-2 & 0 & 1 & 0 & 0 & 0 & 2 & 0 \\
0 &  -2& 1 & 0 & 2 & 0 & 0 & 0 \\
1 & 1 & -2& 1 & 0 & 0 & 0 & 0 \\
0 & 0 & 1 &-2 & 0 & 0 & 0 & 0 \\
0 & 2 & 0 & 0 &-2 & 1 & 0 & 0 \\
0 & 0 & 0 & 0 & 1 &-2 & 0 & 0 \\ 
2 & 0 & 0 & 0 & 0 & 0 &-2 & 1 \\
0 & 0 & 0 & 0 & 0 & 0 & 1 & -2
\end{array}\right)\ {\rm with}\ \det(I)=36.
\end{align*}
Since the determinants are all non-zero, all the cases are impossible.
\qed

\begin{lemm}\label{I1I1} The case (2b) (of two smooth curves intersecting tangentially at two smooth points) in Lemma \ref{2Inv} does not occur.
\end{lemm}
\ni{\bf Proof.} Say $(\Sigma,\Sigma')$ is of type $(I_1,I_1)$ (resp. $(I_1,I_2)$) of (2b) with $\Sigma\cap\Sigma'=\{a,b\}$. Assume that $\Sigma_Z$ intersects $E_{a1}$ and $E_{b1}$. By Lemma \ref{invtor2}, {we can consider the subsystem $\Lambda:=\left<S:=12\Sigma_Z,T:=12\Sigma_Z'\right>\subseteq|24\mu^*L_Y-8E_{a1}-4E_{a2}-8E_{b1}-4E_{b2}|$ (resp. $S:=12\Sigma_Z$, $T=12(\Sigma_Z'+E_{c1}+E_{c2})$) with $S\cdot T=0$.  It defines a morphism $\varphi:Z\rightarrow P^1_\bC$ so that $S$ and $T$ are two special fibers. Since $S\cdot E_{a1}=12$, $\varphi|_{E_{a1}}$ is a degree 12 ramified cover over $P^1_\bC$ with ramification index\footnote{For a map $z\mapsto z^r$, the ramification index is $r-1$.}  11 along $S\cap E_{a1}$ and $T\cap E_{a1}$.  The connected curve $E_{a2}$ is disjoint from $S$ and $T$ and hence sits in the (scheme theoretic) fiber $F:=\varphi^*(\varphi(E_{a2}))$ of $\varphi:Z\rightarrow P^1_\bC$. Since $E_{a2}\cdot E_{a1}=1$ and $F\cdot E_{a1}=12$, either there are more than one components of $F$ passing through $Q:=E_{a2}\cap E_{a1}$ or $12E_{a1}\leq F$. In either cases, the ramification index of $\varphi|_{E_{a1}}$ at $Q$ is at least 1. The contribution of the ramification indices at $S\cap E_{a1},  T\cap E_{a1}$ and $Q$  violates Riemann-Hurwitz formula: $\deg(R_{\varphi|_{E_{a1}}})=12(2)-2=22$. 
\qed

\begin{lemm}\label{NN} The case (3a) (of two double curves intersecting at exactly one point of multiplicity four) in Lemma \ref{2Inv} does not occur.
\end{lemm}
\ni{\bf Proof.} Say $\Sigma=2C_1\sim 2L$ and $\Sigma'=2C_2\sim 2L'$. Then $C_i$'s are $\Aut(M)$-invariant and by Lemma \ref{invtor2} we can consider the subsystem $\Lambda:=\left<24C_1, 24C_2,12(C_1+C_2)\right>\subseteq|24L|$. The same argument as in the proof of Lemma \ref{X} then leads to a contradiction.
\qed

\ms

\ni{\bf 6.6} At this point, we are not able to rule out all possible pairs of invariant curves in Lemma \ref{2Inv}. On the other hand, when $\HH_1(M/\Aut(M),\bZ)$ contains a non-trivial 3-torsion element, there are three cubic roots of $K_M$. Hence we are lead to study triples of three invariant curves. The hypothesis in Lemma \ref{ECeasy} is fulfilled once we rule out all possible triples raised in this way. In the end of this section, we prove the nonexistence of some type of triples under the condition $\Aut(M)=C_3$. For other possibilities, we need more assumptions, cf. Section \ref{Sec7}. 

Suppose that there are three distinct $\Aut(M)$-invariant cubic root $L,L',L''$ of $K_M$ and let $\Sigma\sim 2L, \Sigma'\sim 2L', \Sigma''\sim 2L''$ be three distinct $\Aut(M)$-invariant curves. 
\begin{lemm}\label{I3I3} The case (1b-2) in Lemma \ref{2Inv} does not occur in a triple $(\Sigma,\Sigma',\Sigma'')$.
\end{lemm}
\ni{\bf Proof.} Say $(\Sigma,\Sigma')$ is of type $(I_3,I_3)$ of (1b-2) with $\Sigma\cap\Sigma'=\{a,b\}$ such that $b\in\Sigma$ and $a\in\Sigma'$ are tacnodes. Assume that $\Sigma_Z$ intersects $E_{a1}$ and $E_{b1}$. By Lemma \ref{invtor2}, we find the subsystem 
\begin{eqnarray*}
\Lambda&:=&\left<S:=12(\Sigma_Z+E_{b1}),
T:=12(\Sigma_Z'+E_{a2})\right>\\
&\subseteq&|24\mu^*L_Y-8E_{a1}-4E_{a2}-4E_{b1}-8E_{b2}|
\end{eqnarray*}
 with $S\cdot T=0$}. By the same argument as in the end of the proof of Lemma \ref{I1I1}, this defines a morphism $\varphi:Z\rightarrow P^1_\bC$ such that $\varphi|_{E_{a1}}$ is a degree 12 ramified cover over $P^1_\bC$ with ramification index 11 along $S\cap E_{a1}$ and $T\cap E_{a1}$.

Now consider the types of $(\Sigma,\Sigma'')$ in Lemma \ref{2Inv}:
\begin{enumerate}

\item[(1b-2)] $(I_3,I_3)$: Then $a\in\Sigma''$ is a tacnode, but there is no such local intersection for $a\in\Sigma'\cap\Sigma''$ from Lemma \ref{localint}.

\item[(1c)] $(I_3,I_1)$: The curve $\Sigma''$ has a unique smooth branch tangential to the tacnode $b\in\Sigma$ and $\Sigma\cap\Sigma''=\{a,c\}$. But then $\Sigma'\cap\Sigma''=\{b\}$ and they intersect transversally, which violates $\Sigma'\cdot\Sigma''=4.$ 

\item[(1c)] $(I_3,I_3)$: The local picture at $b\in\Sigma\cap\Sigma'\cap\Sigma''$ is the same as in the case of $(I_3,I_1)$ in  (1c). Hence the same argument as above leads to a contradiction.
\end{enumerate}
The only possibility left is when $\Sigma''$ is of type $(N)$ so that $\Sigma''$ intersects $\Sigma$ transversally at $b$ and tangentially to $\Sigma'$ at $b$, or $\Sigma''$ intersects $\Sigma$ tangentially at $a$ and transversally to $\Sigma'$ at $a$. Assume that we are in the latter case and $\Sigma''_Z=2C''_Z$. Then the connected curve $C''_Z$ is disjoint from $S$ and $T$ and hence sits in the (scheme theoretic) fiber $F:=\varphi^*(\varphi(C''_Z))$ of $\varphi:Z\rightarrow\bP^1$. Since $C''_Z\cdot E_{a1}=1$ and $F\cdot E_{a1}=12$, either there are more than one components of $F$ passing through  $Q:=C''_Z\cap E_{a1}$ or $12C''_Z\leq F$. In either cases, the ramification index of $\varphi|_{E_{a1}}$ at $Q$ is at least 1. The same argument as in Lemma \ref{I1I1} leads to a contradiction to Riemann-Hurwitz formula. The other case is treated similarly.	
\qed

\begin{lemm}\label{1b3} The case (1b-3) in Lemma \ref{2Inv} does not occur in a triple $(\Sigma,\Sigma',\Sigma'')$.
\end{lemm}
\ni{\bf Proof.} Say $(\Sigma,\Sigma')$ is of type $(I_1,I_3)$ of (1b-3) with $\Sigma\cap\Sigma'=\{a,b\}$ such that $b\in\Sigma'$ is a tacnode. Assume that $\Sigma_Z$ intersects $E_{a1}$ and $E_{b1}$. By Lemma \ref{invtor2}, we find the divisors $S:=12\Sigma_Z,T:=12(\Sigma_Z'+E_{b2})\in|24\mu^*L_Y-8E_{a1}-4E_{a2}-8E_{b1}-4E_{b2}|$ with $S\cdot T=0$. The pencil $\Lambda:=\left<S,T\right>$ defines a morphism $\varphi:Z\rightarrow P^1_\bC$ such that $\varphi|_{E_{a1}}$ is a degree 12 ramified cover over $P^1_\bC$ with ramification index 11 along $S\cap E_{a1}$ and $T\cap E_{a1}$. The connected curve $E_{a2}$ is disjoint from $S$ and $T$ and hence sits in the (scheme theoretic) fiber $F:=\varphi^*(\varphi(E_{a2}))$ of $\varphi:Z\rightarrow P^1_\bC$. Since $E_{a2}\cdot E_{a1}=1$ and $F\cdot E_{a1}=12$, either there are more than one components of $F$ passing through $Q:=E_{a2}\cap E_{a1}$ or $12E_{a2}\leq F$. In either cases, the ramification index of $\varphi|_{E_{a1}}$ at $Q$ is at least 1. Again this violates Riemann-Hurwitz formula as explained in the proof of Lemma \ref{I1I1}. 

Suppose now $(\Sigma,\Sigma')$ is of type $(I_2,I_3)$ of (1b-3) with $\Sigma\cap\Sigma'=\{a,b\}$ such that $b\in\Sigma'$ is the tacnode and $c\in\Sigma$ is the node. Assume that $\Sigma_Z$ intersects $E_{a1}$ and $E_{b1}$. By Lemma \ref{invtor2}, we find the divisors $S:=12(\Sigma_Z+E_{c1}+E_{c2}),T:=12(\Sigma_Z'+E_{b2})\in|24\mu^*L_Y-8E_{a1}-4E_{a2}-8E_{b1}-4E_{b2}|$ with $S\cdot T=0$. Similarly the pencil $\Lambda:=\left<S,T\right>$ defines a morphism $\varphi:Z\rightarrow P^1_\bC$ such that $\varphi|_{E_{a1}}$ is a degree 12 ramified cover over $P^1_\bC$ with ramification index 11 along $S\cap E_{a1}$ and $T\cap E_{a1}$. Consider the connected curve $E_{a2}$ and the (scheme theoretic) fiber $F:=\varphi^*(\varphi(E_{a2}))$ of $\varphi:Z\rightarrow P^1_\bC$. The same argument as in the proof of Lemma \ref{I1I1} leads to a contradiction to Riemann-Hurwitz formula.
\qed

\bs

\section{The case of ${\rm Aut}(M)=C_3$ and $\HH_1(M/C_3,\bZ)=C_3$}\label{Sec7}

To illustrate how Section \ref{Sec6} helps us to prove the Main Theorem, we focus on the case when a fake projective plane $M$ has ${\rm Aut}(M)=C_3$ and $\HH_1(Y,\bZ)=C_3$, where $Y=M/\Aut(M)$.  Since there are three distinct $\Aut(M)$-invariant cubic root $L,L',L''$ of $K_M$ from the discussion in Section \ref{Sec2}, by Lemma \ref{ECeasy} it is enough to show that one of $2L, 2L', 2L''$ has no global sections. We assume the contrary and let $\Sigma\sim 2L, \Sigma'\sim 2L', \Sigma''\sim 2L''$ be three distinct $\Aut(M)$-invariant curves. 

\begin{lemm}\label{NX} If $M$ is a fake projective plane with $\Aut(M)=C_3$ and $\HH_1(Y,\bZ)=C_3$, then an invariant curve $\Sigma\sim 2L$ of type $(N)$ or $(X)$ does not exist.
\end{lemm}
\ni{\bf Proof.} If $\Sigma$ is of type $(N)$ and $\Sigma=2C\equiv 2L$, then $C$ is invariant and $C\sim L+\alpha$ for some $3$-torsion $\alpha$ coming from $\HH_1(Y,\bZ)=C_3.$ Hence $3\alpha$ can only be trivial and $h^0(M,L+\alpha)\neq0$. This contradicts to $p_g(M)=h^0(M,3L)=h^0(M,3(L+\alpha))=0.$ Similarly, if $\Sigma=\Sigma_1+\Sigma_2$ is of type $(X)$, then $\Sigma_i$ being invariant must be in the class $L+\omega_i$, where $\omega_i$ is an invariant 3-torsion. But then $3\Sigma_i\sim 3L=K_M$, contradicting to $p_g(M)=0.$
\qed

\begin{rem} Compare to the proof of Lemma \ref{NX}, the proof of Lemma \ref{X} shows that in general it is harder to rule out all possible invariant curves in the numerical class $2L$.
\end{rem}

As a consequence of Section \ref{Sec6} and Lemma \ref{NX}, to consider a triple of three distinct invariant curves of numerical type $2L$, we only need to consider pairs of invariant curves of the following types in Lemma \ref{2Inv}:
	\begin{enumerate}		
		\item[(1a)] $(I_1,I_2)$ with intersection configuration $3(tan-node)+1(tr)$;
		\item[(1c)] $(I_1,I_3)$, $(I_3,I_3)$ with intersection configuration $4(tan-tac)$;
	\end{enumerate}
	
\begin{lemm}\label{I1I21a} Type $(I_1,I_2)$ of (1a) in Lemma \ref{2Inv} does not occur in a triple $(\Sigma,\Sigma',\Sigma'')$.
\end{lemm}
\ni{\bf Proof.} If $(\Sigma,\Sigma')$ is of type $(I_1,I_2)$ of (1a) and say $a\in\Sigma'$ the node, then $(\Sigma,\Sigma'')$ can only be of type $(I_1,I_2)$ of (1a) or $(I_1,I_3)$ of (1c).  In the former case, $(\Sigma',\Sigma'')$ is of type $(I_2,I_2)$ and must have $\Sigma'\cdot\Sigma''>4$, which is absurd. In the latter case, we must have that $b\in\Sigma''$ is the tacnode. But then it is only possible that $(\Sigma',\Sigma'')$ to be of type $(I_2,I_3)$ in (1b-3) with $b=(tr-tac)$ and $c=(tan-sm)$ for the third fixed point $c$ of $\Aut(M)$. This case has been ruled out in Lemma \ref{1b3}.
\qed

\ms

Up to this point, we see that the triple $(\Sigma, \Sigma',\Sigma'')$, up to reordering, can only be of type $(I_1,I_3,I_3)$ or $ (I_3,I_3,I_3),$ where $(I_1,I_3)$ and $(I_3,I_3)$ are both from $(1c)$.

\begin{lemm} The type $(I_1,I_3,I_3)$ does not occur.
\end{lemm}
\ni{\bf Proof.} Say we have type $(I_1,I_3,I_3)$ such that $b\in\Sigma'$ and $c\in\Sigma''$ are tacnodes and $\{b,c\}=\Sigma'\cap\Sigma''$. But both $(\Sigma,\Sigma')$ and $(\Sigma,\Sigma'')$ being type $(I_1,I_3)$ in (1c) then implies that $\Sigma$ is tangential to $\Sigma''$ along $b$ and $c$, which leads to
$$4=\Sigma\cdot\Sigma''=\mult_b(\Sigma\cap\Sigma'')+\mult_c(\Sigma\cap\Sigma'')\geq1+4.$$
\qed

\begin{lemm}\label{3I3} The type $ (I_3,I_3,I_3)$ does not occur.
\end{lemm}
\ni{\bf Proof.} Say $\Sigma=(s_1=0)$, $\Sigma'=(s_2=0)$, and $\Sigma''=(s_3=0)$ on $M$. Here $s_i$ descends to $\ts_i\in H^0(Y,2L_i)$ and pulls back to $t_i:=\mu^*\ts_i\in H^0(Z,\mu^*(2L_i))$, for $i=1,2,3$.  As all invariant torsions has order three from the proof of Lemma \ref{invtor}, we can form the linear system $\Lambda:=\left<t_1^3,t_2^3,t_3^3,t_1t_2t_3\right>$ on $Z$. Note that by construction $\Lambda\subseteq|2K_Z|.$

It is easy to see that these four sections are linearly independent: Suppose this is not true. Clearly there is an induced relation $A s_1^3+B s_2^3+C s_3^3+Ds_1s_2s_3=0$ for some constants $A,B,C,D$. Then considering a generic point of $s_1=0$, $Bs_2^3+Cs_3^3=0$ along $s_1=0$.  But $s_2, s_3$ have different vanishing orders at the fixed points.  Hence $B=C=0$ (say, by taking derivatives along the curve $s_1=0$). Now we have $As_1^3+Ds_1s_2s_3=0$, or $As_1^2+Ds_2s_3=0$.  But again at a generic point of $s_1=0$, $s_2$ and $s_3$ do not vanish. This forces $D=0$.  Hence $A=0$ as well. This also shows that $\Lambda=|2K_Z|$ as $h^0(Z,2K_Z)=4$ by Riemann-Roch formula and Kawamata-Viehweg vanishing.

Note that $\Lambda$ is base point free as $\Sigma_Y\cap\Sigma'_Y\cap\Sigma''_Y=\emptyset$ and defines a morphism $\Phi:Z\rightarrow P^3_\bC$ via 
$$[x:y:z:w]=[t_1^3:t_2^3:t_3^3:t_1t_2t_3].$$
Clearly, $\Phi(Z)=S=(w^3=xyz)$, which is a singular normal cubic surface with three $A_2$ singularities. Since $K_Z\cdot E=0$ exactly along $\mu$-exceptional curves, there is a factorization $\Phi:Z\xrightarrow{\mu} Y\xrightarrow{\Psi} S$. As $2K_Z=\Phi^*\cO_S(1)$, $\Psi$ is finite with 
$$\deg(\Psi)=\frac{(2K_Z)^2}{\cO_S(1)^2}=4.$$
Note that $Y$ and $S$ are Gorenstein and being a cubic $K_S=\cO_S(-1)$. By Riemann-Hurwitz formula, $K_Y=\Psi^*K_S+R$. On one hand, local computation along generic point of $(t_i=0)$ gives the ramification index two. Hence 
$$R\geq D:=2((t_1=0)+(t_2=0)+(t_3=0))\equiv4K_Y.$$
On the other hand, as $\Psi^*K_S=-\Psi^*\cO_S(1)=-2K_Y$, we get 
$$3=K_Y^2=K_Y\cdot (\Psi^*K_S+R)=K_Y\cdot (2K_Y+(R-D))=6+K_Y\cdot(R-D)\geq6.$$
This is a contradiction.
\qed

\ms

We are now ready to prove the Main Theorem for fake projective planes in line 7 to 9 of Table \ref{Tab:Main}.

\begin{coro}\label{C3a} Suppose that $M$ is a fake projective plane with $\Aut(M)=C_3$ and $\HH_1(M/\Aut(M),\bZ)= C_3$. There is an invariant cubic root $L$ with $K_M=3L$ such that the sequence $\cO_M,-L,-2L$ forms an exceptional collection.
\end{coro}

\ni{\bf Proof.} By Lemma \ref{invcubicroot} and \ref{invtor}, there are three $\Aut(M)$-invariant cubic roots $L,L',L''$ of $K_M$, i.e., $K_M=3L=3L'=3L''.$ 
If all of $2L,2L'$,and $2L''$ have non-trivial global sections, then by Lemma \ref{invcurve1} there is a triple $(\Sigma,\Sigma',\Sigma'')$ of distinct $\Aut(M)$-invariant curves in the numerical class $2L$. However, from results in Section \ref{Sec6} and Lemma \ref{NX} to Lemma \ref{3I3} in this section, all possible configurations of such a triple are ruled out and this is absurd. Hence one of $2L, 2L'$, or $2L''$ has no global section. It follows that from Lemma \ref{ECeasy} there exists an exceptional collection of the expected type. 
\qed
\bs

\section{The case of ${\rm Aut}(M)=C_3$ and $\HH_1(M/C_3,\bZ)=C_2\times C_3$}\label{Sec8}

Suppose that $M$ is a fake projective plane with ${\rm Aut}(M)=C_3$ and $\HH_1(M/\Aut(M),\bZ)=C_2\times C_3$.  There are 6 classes of such $M$ in Table \ref{Tab:Main} and 3 classes in Table \ref{Tab:unproven}. First of all, by using the same trick as in the last section, we show that there cannot be triples of different invariant curves in numerical class $2L$ except one case, cf Proposition \ref{last}. This implies the vanishing of many global sections for invariant line bundles of numerical type $L$ or $2L$. If $K_M=3L$, then we find an exceptional collection, cf. Corollary \ref{C3b}. For fake projective planes not treated in Corollary \ref{C3a} and \ref{C3b}, we refer the reader to Section \ref{Sec9}.

\ms

\ni{\bf 8.1} Again, we will consider a triple of three distinct $\Aut(M)$-invariant curves $(\Sigma, \Sigma',\Sigma'')$ in the numerical class $2L$. Note that Lemma \ref{NX} does not apply here and hence we cannot quote results in Section \ref{Sec7} directly. From Section \ref{Sec6}, the intersection type of any two invariant curves can only be from the following list: 
	\begin{enumerate}		
		\item[(1a)] $(I_1,I_2)$ with intersection configuration $3(tan-node)+1(tr)$;
		\item[(1b-1)] $(N,I_3)$ with intersection configuration $4=2\cdot2(tr-tac)$;
		\item[(1c)] $(I_1,I_3)$, $(I_3,I_3)$ with intersection configuration $4=4(tan-tac)$;
		\item[(3b)] $(N,I_1)$, $(N,I_2)$ with intersection configuration $4=2\cdot(1(tr)+1(tr))$;
		\item[(3c)] $(N,I_1)$, $(N,I_3)$ with intersection configuration $4=2\cdot2(tan-sm)$.
	\end{enumerate}

\begin{prop}\label{last} Suppose that $M$ is a fake projective plane with $\Aut(M)=C_3$ and $\HH_1(M/\Aut(M),\bZ)=C_2\times C_3$. Then there exists no triple of three distinct $\Aut(M)$-invariant curves $(\Sigma, \Sigma',\Sigma'')$ in the numerical class $2L$ except, up to reordering, when the type is $(N, I_1, I_2)$:
\begin{center}
\includegraphics[scale=0.7]{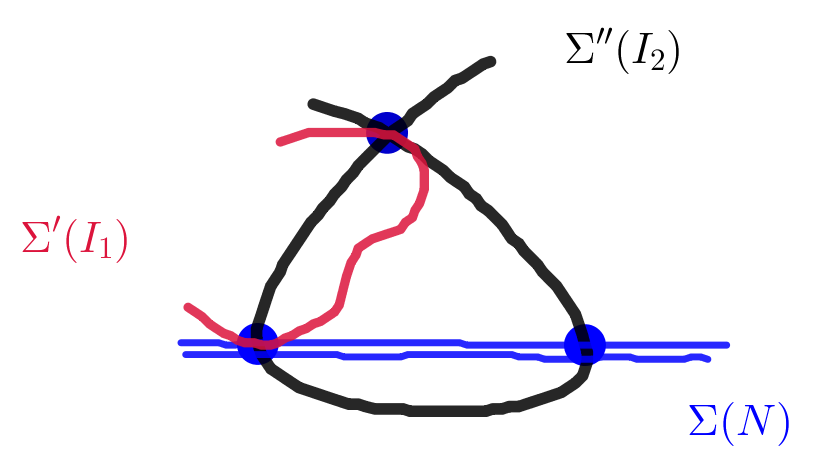}
\end{center}
\end{prop}
\ni{\bf Proof.} As discussed above, we only need to consider the case when exactly one of these invariant curves is of type $(N)$, say $\Sigma$. By Lemma \ref{NN}, none of $\Sigma'$ and $\Sigma''$ can be of type $(N)$. Hence from the above list, 
the type of $(\Sigma',\Sigma'')$ can only be (1a) or (1c). If $(\Sigma',\Sigma'')$ is of the type $(I_1,I_2)$ in (1a), then we can only have a triple of type $(N,I_1,I_2)$ as depicted in the statement. Hence hereafter we assume that $(\Sigma',\Sigma'')$has intersection configuration $4=4(tan-tac)$ from (1c) of Lemma \ref{2Inv}. 

Assume that $b\in\Sigma''$ is the tacnode and $\Sigma'$ intersects $\Sigma''$ transversally at $b$ so that $\Sigma'\cap\Sigma''=\{b\}$. Suppose that $\{a,b\}\subseteq\Sigma''$ and $\{b,c\}\subseteq\Sigma'$. Note that $\Sigma$ does not pass through $a,b$ simultaneously since there is no such type $(N,I_3)$ of $(\Sigma, \Sigma'')$ in the above list.

If $\Sigma$ passes through $a, c$, then $\Sigma\cap\Sigma'\cap\Sigma''=\emptyset$. In particular, a similar argument as the proof of Lemma \ref{3I3} applies: Since $6\tau=0$ for all $\Aut(M)$-invariant torsion line bundles, there is a linear system $\Lambda:=\left<t_1^6,t_2^6,t_3^6,t_1^2t_2^2t_3^2\right>\subseteq|4K_Z|.$ The same computation as in the proof of Lemma \ref{3I3} leads to a contradiction from $$K_Y=\Psi^*K_S+R=-4K_Y+(R-D)+D,$$
where $$D:=5((t_1=0)+(t_2=0)+(t_3=0))\equiv 10K_Y.$$
Hence it is only possible that $\{b,c\}\subseteq\Sigma$ and $(\Sigma',\Sigma'')$ is of type $(I_1,I_3)$ in $(3b)$. 
In particular, $\Sigma\cap\Sigma'=\{b,c\}.$
	
Assume that  $\Sigma_Z$ intersects $E_{b1}$ and $E_{c1}$. Since $6\tau=0$ for all $\Aut(M)$-invariant torsion line bundles,  we can consider the subsystem $$\left<S:=6(\Sigma_Z+E_{b1}+E_{c1}), T:=6\Sigma_Z'\right>\subseteq|12\mu^*L_Y-2E_{b1}-4E_{b2}-2E_{c1}-4E_{c2}|,$$ where $S\cdot T=0$. As in the end of the proof of Lemma \ref{I1I1}, this defines a morphism $\varphi:Z\rightarrow P^1_\bC$ such that $\varphi|_{E_{b2}}$ is a degree 6 ramified cover over $P^1_\bC$ with ramification index 5 along $S\cap E_{b2}$ and $T\cap E_{b2}$.  The connected curve $\Sigma''_Z$ is disjoint from $T$ and hence sits in the (scheme theoretic) fiber $F:=\varphi^*(\varphi(\Sigma''_Z))$ of $\varphi:Z\rightarrow P^1_\bC$. Since $\Sigma''_Z\cdot E_{b2}=2$ and $F\cdot E_{b2}=6$, either there are more than one components of $F$ passing through one of $\{Q,Q'\}:=\Sigma''_Z\cap E_{b2}$ or $3\Sigma''_Z\leq F$. In either cases, the ramification index of $\varphi|_{E_{b2}}$ at one of $\{Q,Q'\}$ is at least 1. The count of the ramification indices at the three points above violates Riemann-Hurwitz formula as in the proof of
Lemma \ref{I1I1}.\qed

\ms

\ni{\bf 8.2} Write $\HH_1(M/\Aut(M),\bZ)=\left<\tau,\omega\right>$, where $\left<\tau\right>\cong C_2$ and $\left<\omega\right>\cong C_3$. We will identify $\tau=(1,0)$, $\omega=(0,1)$, and use the additive notion. Write $K_M=3L+\mu$, where $L$ is a fixed $\Aut(M)$-invariant line bundle and $\mu$ is some $\Aut(M)$-invariant torsion. From the proof of Lemma \ref{invtor} and by abuse of notion, we can assume that $\mu\in\left<\tau,\omega\right>.$
There are two cases:
\begin{enumerate}
	\item $M$ is not in the class $\cC_{18}$: By Lemma \ref{cubicroot}, we can put $\mu=0.$
	\item $M$ is in the class $\cC_{18}$: As $3L+\tau=3(L+\tau)$, we can always choose $\omega$ corresponding to a generator of the $C_3$-factor of $\HH_1(M/\Aut(M),\bZ)$ so that $\mu=\omega$.
\end{enumerate}
 Hence hereafter we fix the setup: $$\HH_1(M/\Aut(M),\bZ)=\left<\tau,\omega\right>\ {\rm and}\ K_M=3L+\omega,$$
 where $\omega=0$ if $M$ is not in the class $\cC_{18}.$

\begin{lemm}\label{van3tor} One of $2L, 2(L+\omega)$, or $2(L+2\omega)$ has no global sections. 
\end{lemm}
\ni{\bf Proof.} Suppose that the contrary holds. From Proposition \ref{last}, there is a triple of invariant curves of type $(N,I_1,I_2)$ with $\Sigma=2C\sim 2(L+k\omega)$ being of type $(N)$ for some $k\in\{0,1,2\}$.  Rewrite $K_M=3(L+k\omega)+\omega$, we may assume that $k=0$ and $\Sigma=2C\sim 2L$. From Table \ref{intnumZ}, $C_Z$ is a smooth elliptic curve and $3\Sigma'\sim 3\Sigma''\sim 6L$ restricts to two sections $t_1:=\mu^*(3\Sigma'_Y)|_{C_Z}, t_2:=\mu^*(3\Sigma''_Y)|_{C_Z}\in H^0(C_Z,\mu^*6L_Y|_{C_Z}).$ From Riemann-Roch formula, it is easy to see that $h^0(C_Z,\mu^*6L_Y|_{C_Z})=2$ and is generated by global sections. On the other hand, $t_1$ and $t_2$ are linearly independent from the description of these three curves in Proposition \ref{last} and hence $|\mu^*6L_Y|_{C_Z}|=\left<t_1,t_2\right>$. However, $\left<t_1,t_2\right>$ has base points along $\mu^{-1}(C_Y\cap\Sigma'_Y\cap\Sigma''_Y)$. This is a contradiction.
\qed

\ms

We can now prove the Main Theorem for fake projective planes in the last six lines of Table \ref{Tab:Main}.

\begin{coro}\label{C3b} Suppose that $M$ is a fake projective plane with $\Aut(M)=C_3$ and $\HH_1(M/\Aut(M),\bZ)=C_2\times C_3$. If $M$ is not in class $\cC_{18}$, then there is an $\Aut(M)$-invariant line bundle $L$ with $K_M=3L$ such that the sequence $\cO_M,-L,-2L$ forms an exceptional collection.
\end{coro}
\ni{\bf Proof.} We have $K_M=3L$ for some $\Aut(M)$-invariant line bundle from Lemma \ref{invcubicroot}, and $L, L+\omega, L+2\omega$ are three distinct invariant cubic roots of $K_M.$ 
Hence the corollary follows immediately from Lemma \ref{van3tor} and Lemma \ref{ECeasy}.
\qed

\bs

\section{Remarks on the other cases with ${\rm Aut}(M)=C_3$}\label{Sec9}
The list of all the fake projective planes with a non-trivial automorphism which are not treated in the Main Theorem is in Table \ref{Tab:unproven}. There are in total 36 non-biholomorphic 
of such fake projective planes. With minor modifications, the results in Section \ref{Sec6} to \ref{Sec8} for a pair or a triple of different invariant curves apply to \emph{all} fake projective planes with non-trivial automorphisms, except  Lemma \ref{NX} and Lemma \ref{van3tor}. The difficulty to prove Conjecture \ref{numGKMS} in general is to establish the following two key ingredients as done in Corollary \ref{C3b}.

The first ingredient is very technical and is the main difficulty in applying our approach to the remaining fake projective planes: Lemma \ref{van3tor} holds if $3\Sigma'_Y\sim3\Sigma''_Y$, or slightly weaker $\mu^*(3\Sigma'_Y)|_{C_Z}\sim\mu^*(3\Sigma''_Y)|_{C_Z}$. This is applicable to fake projective planes whose $\HH_1(M,\bZ)$ contains a  unique $C_3$-factor, for which we consider three $\Aut(M)$-invariant cubic roots of $K_M$ as in Lemma \ref{van3tor}. However, the proof does not work for a general choice of three invariant line bundles. In particular, Lemma \ref{van3tor} fails for any choice of three invariant line bundles when $\HH_1(M/\Aut(M),\bZ)$ contains no $C_3$-factor. To tackle the first difficulty, we propose the following question. 

\begin{ques}\label{q} Let $M$ be a fake projective plane with $\Aut(M)=C_3$ and $K_M\equiv 3L$ for an $\Aut(M)$-invariant line bundle $L$.  Is it true that there exists no triple $(\Sigma,\Sigma',\Sigma'')$ of distinct $\Aut(M)$-invariant curves of type $(N,I_1,I_2)$ in the numerical class $2L$ as in Proposition \ref{last}?
\end{ques}

A positive answer to Question \ref{q} does not prove Conjecture \ref{numGKMS} directly, but shall be taken as a weak solution to it.

The second ingredient is that, even if Lemma \ref{van3tor} holds for the choice of three invariant cubic roots of $K_M$, we still need $K_M=3L$ to apply Lemma \ref{ECeasy}. This is the main difficulty to prove Conjecture \ref{numGKMS} for fake projective planes in class $\cC_{18}$. 
If $M$ possesses many invariant torsions, then we may apply our approach with the following generalization of Lemma \ref{ECeasy}.

\begin{lemm}\label{EC} For a choice of torsions $\mu_1,\mu_2, \omega$ on a fake projective plane $M$ such that $K_M=3L+\omega$, the sequence $\cO_M,-(L+\mu_1), -(2L+\mu_2)$ forms an exceptional collection if and only if 
$$h^0(M, L+\mu_1)=h^0(M, L+\omega-\mu_2)=h^0(M, L+\mu_2-\mu_1)=0$$
and
$$h^0(M, 2L+\omega-\mu_1)=h^0(M, 2L+\mu_2)=h^0(M, 2L+\omega+\mu_1-\mu_2)=0.$$
\end{lemm}
\ni{\bf Proof.} The required vanishing for the given sequence of line bundles to be an exceptional collection is given by 
$$	\begin{cases} h^0(M, L+\mu_1)=h^1(M, L+\mu_1)=h^2(M, L+\mu_1)=0,\\
	h^0(M, 2L+\mu_2)=h^1(M, 2L+\mu_2)=h^2(M, 2L+\mu_2)=0,\ {\rm and}\\
	h^0(M, L+\mu_2-\mu_1)=h^1(M, L+\mu_2-\mu_1)=h^2(M, L+\mu_2-\mu_1)=0.\\
	\end{cases}$$
	By Serre duality, e.g. $h^2(M, L+\mu_1)=h^0(M, 2L+\omega-\mu_1)$, this gives the necessary condition. Conversely, together with $p_g(M)=q(M)=0$ and $\chi(L')=\chi(2L')=1$ for any positive line bundle $L'$ generating $\NS(M)_\bQ$, the prescribed vanishing of $h^0$ implies all the required vanishing of $h^1$.
\qed

\ms

To apply Lemma \ref{EC}, one can consider all possible invariant numerical torsions $(\mu_1,\mu_2,\omega)$ and apply the discussion in Section \ref{Sec6} to \ref{Sec8} to obtain some vanishing of invariant global sections. 
However, we have checked that in some cases there is no compatible choice of torsion line bundles to generate all the required vanishing conditions.

{\small
\centering
\begin{table}
   \captionof{table}{FPP with $\Aut\neq\{1\}$ not covered in Main Theorem\label{Tab:unproven}}
\begin{adjustbox}{center}$
\begin{array}{|c|c|c|c|c|c|}
\hline
\mbox{class}&M&{\rm Aut}(M)&\HH_1(M,\bZ)&H& \HH_1(M/H,\bZ)\\ \hline\hline
(a=1,p=5,\emptyset)&(a=1,p=5,\emptyset, D_3)&C_3&C_2\times C_4\times C_{31}&C_3&C_2\times C_4\\ \hline
(a=1,p=5,\{2\})&(a=1,p=5,\emptyset, \{2\},D_3)&C_3&C_4\times C_{31}&C_3&C_4\\ \hline
(a=2,p=3,\{2\})&(a=2,p=3,\{2\}, D_3)&C_3&C_2^2\times C_{13}&C_3&C_2\times C_2\\ \hline
(a=2,p=3,\emptyset)&(a=2,p=3,\emptyset, D_3)&C_3&C_2^2\times C_{13}&C_3&C_2\times C_2\\ \hline
(a=7,p=2,\emptyset) &(a=7,p=2,\emptyset, D_3 X_7)&C_3&C_2\times C_7&C_3&C_2\\ \hline
(a=7,p=2,\{7\}) &(a=7,p=2,\{7\}, D_3 7_7)&C_3&C_2\times C_7&C_3&C_2\\ \cline{2-6}
&(a=7,p=2,\{7\}, D_3 7'_7)&C_3&C_2^2\times C_7&C_3&C_2\times C_2\\ \hline
(a=7,p=2,\{3\})&(a=7,p=2,\{3\},D_3)&C_3&C_2\times C_4\times C_7&C_3&C_2\times C_4\\ \hline
(a=7,p=2,\{3,7\})&(a=7,p=2,\{3,7\},D_3)&C_3&C_4\times C_7&C_3&C_4\\ \hline
(a=15,p=2,\emptyset)&(a=15,p=2,\emptyset,D_3)&C_3&C_2^2\times C_7&C_3&C_2\times C_2\\ \hline
(a=15,p=2,\{5\})&(a=15,p=2,\{5\},D_3)&C_3&C_2\times C_7&C_3&C_2\\ \hline
(\cC_{10},p=2,\emptyset)&(\cC_{10},p=2,\emptyset, D_3)&C_3&C_2\times C_7&C_3&C_2\\  \hline
(\cC_{10},p=2,\{17-\})&(\cC_{10},p=2,\{17-\}, D_3)&C_3&C_7&C_3&\{1\}\\  \hline
(\cC_{18},p=3,\{2\})&(\cC_{18},p=3,\{2\},D_3)&C_3&C_2\times C_3\times C_{13}&C_3&C_2\times C_3\\  \cline{2-6}
&(\cC_{18},p=3,\{2\},(dD)_3)&C_3&C_2\times C_3&C_3&C_2\times C_3\\  \cline{2-6}
&(\cC_{18},p=3,\{2\},(d^2D)_3)&C_3&C_2\times C_3&C_3&C_2\times C_3\\  \hline
(\cC_{20},\{v_2\},\{3+\})&(\cC_{20},\{v_2\},\{3+\},D_3)&C_3&C_4\times C_7&C_3&C_4\\ \hline
(\cC_{20},\{v_2\},\{3-\})&(\cC_{20},\{v_2\},\{3-\},D_3)&C_3&C_4\times C_7&C_3&C_4\\ \hline
\end{array}$
\end{adjustbox}
\end{table}
}

\bs

\ni{\bf Acknowledgment.} The authors would like to thank the referees for helpful comments. This work is partially done during the first author's visit at National Center of Theoretical Science in Taiwan, and the second author's visit of the Institute of Mathematics of the University of Hong Kong during the summers of 2013 and 2014. The authors thank the warm hospitality of the institutes. The first author is supported by the grant MOST 108-2115-M-006-016 and an internal grant of National Cheng Kung University. The second author was partially supported by a grant from the National Science Foundation.

\bibliographystyle{plain}
\bibliography{FPP.bib}
\end{document}